\documentclass{bcp}
\usepackage{amsmath,amsthm}
\usepackage{amssymb}

\usepackage{graphicx}

\usepackage[T1]{fontenc}

\newtheorem{theorem}{Theorem}[section]
\newtheorem{proposition}[theorem]{Proposition}
\newtheorem{lemma}[theorem]{Lemma}
\newtheorem{corollary}[theorem]{Corollary}

\theoremstyle{definition}
\newtheorem{definition}[theorem]{Definition}
\newtheorem{example}[theorem]{Example}
\newtheorem{remark}[theorem]{Remark}

\DeclareMathOperator{\End}{End}

\newcommand\lie[1]{\mathfrak{#1}}
\newcommand{\g}{\lie{g}}

\renewcommand{\O}{\mathcal{O}}

\newcommand{\G}{\varGamma}

\def \N {{\mathcal {N}}}
\def \d {{\mathrm {d}}}

\begin{document}

\keywords{Poisson-Lie groups, Lie groupoids, Lie algebroids, Poisson groupoids, symplectic groupoids, multiplicative multivectors, mutiplicative forms, generalized geometry, Dirac structures, generalized complex structures.}
\mathclass{Primary 53D17, 22A22; Secondary 70G45.}

\abbrevauthors{Y. Kosmann-Schwarzbach}
\abbrevtitle{Multiplicativity}

\title{Multiplicativity, from Lie groups to generalized geometry}

\author{Yvette Kosmann-Schwarzbach}
\address{Paris, France\\
E-mail: yks@math.cnrs.fr}

\maketitlebcp

\begin{abstract}
We survey the concept of multiplicativity from its initial appearance in the theory of {Poisson-Lie groups} to the far-reaching generalizations, for multivectors and differential forms in the geometry and the generalized geometry of Lie groupoids, as well as their infinitesimal counterparts in the theory of Lie algebroids. 
\end{abstract}

\section{Introduction}

The aim of this survey is to trace the development of the concept of multiplicativity, from its appearance in Lie-group theory in 1983, when Vladimir Drinfeld \cite{D83} defined the Lie groups with a compatible Poisson structure that were later called ``Poisson-Lie groups'', to the theory of Lie groupoids and Lie algebroids and their generalized geometry. On a Lie group, there is a multiplication. It is natural to require that
any additional structure on the group, such as a Poisson structure, be compatible with that multiplication.
{Multiplicativity} is what expresses this compatibility.

In 1987, at the time that Kirill Mackenzie was publishing his first book on Lie groupoids and Lie algebroids, groupoids entered the picture in a very important way as Alan Weinstein \cite{W87}, Mikhail Karasev and V. P. Maslov \cite{K86} \cite{KM93}, 
and  {Stanis{\l}aw Zakrzewski} \cite{Z90}, independently, introduced {\it symplectic groupoids}. There followed Wein\-stein's more general theory of {\it Poisson groupoids} in 1988 \cite{W88}. 
A Poisson groupoid is a Lie groupoid equipped with a multiplicative Poisson bivector, and, in particular, a symplectic groupoid is a Lie groupoid equipped with a multiplicative, closed and non-degenerate $2$-form.
A decisive development was the proof by Mackenzie in 1992 \cite{M92} that a bivector on a Lie group is multiplicative if and only if it defines a morphism from the cotangent groupoid of the Lie group to its tangent group considered as a groupoid over a point.                     
Once the groupoid structure of both the tangent and the cotangent bundles of Lie groups and groupoids were defined, the definition of multiplicative bivectors on Lie groups, suitably formulated as a morphism property, could be extended in three directions, (i) passing from the case of Lie groups to the more general case of Lie groupoids, (ii) passing from the case of bivector fields to that of multivector fields of any degree, and (iii) treating the dual case, that of differential $2$-forms and, more generally, that of differential forms of arbitrary degree. As a rule, we shall write ``bivector'' for ``bivector field'', ``multivector'' for ``multivector field'' and ``form'' for ``differential form''. 

Lie algebroids are the infinitesimal counterparts of Lie groupoids. What are the infinitesimal counterparts that arise in the differentiation process of the multiplicative tensors on Lie groupoids? 
Applying the Lie functor to a multiplicative multivector on a Lie groupoid yields a multivector on its Lie algebroid that not only is ``linear'', i.e., is characterized in terms of vector-bundle morphisms, but also satisfies a property of ``infinitesimal multiplicativity''. Such multivectors have been called ``morphic'' \cite{MX98} \cite{BCO09} because they are characterized in terms of Lie-algebroid morphisms \cite{HM90}. For example, the Lie functor applied to a multiplicative bivector defining a Poisson-Lie structure on a Lie group yields a linear Poisson bivector on its Lie algebra, $\g$, which defines a Lie algebra structure on the dual vector space, $\g^*$. Because this Poisson bivector is not only linear but also morphic, it satisfies a compatibility condition and the pair $(\g,\g^*)$ is a Lie bialgebra.
This argument extends to Lie bialgebroids which are the infinitesimals of Poisson groupoids. 
Similarly the Lie functor applied to multiplicative forms on a Lie groupoid yields ``morphic'' forms on its Lie algebroid. In the case of a morphic, i.e., 
``infinitesimally multiplicative'', $2$-form on a Lie algebroid, $A \to M$, there is an associated vector-bundle morphism from $A$ to the cotangent bundle, $T^*M$, of $M$ which has been called in the literature an ``IM-form''  \cite{BCO09} because, in the case of a closed or relatively closed $2$-form, it determines the morphic 2-form \cite{BCWZ04}.
The theorem on the integration of a morphic multivector (respectively, form) on the Lie algebroid of a source-simply-connected Lie groupoid into a multiplicative multivector (respectively, form) on the groupoid can then be recognized as a particular case of the integration 
of morphisms of Lie algebroids into Lie-groupoid morphisms \cite{MX00}.

We first recall Poisson-Lie groups and the definition and properties of multiplicative multivectors on Lie groups in Section \ref{groups}. We introduce groupoids, their tangent and cotangent bundles in Section \ref{groupoids}, and we describe their  role in the theory of Poisson-Lie groups.
In Section \ref{PG}, we deal with the theory of symplectic and Poisson groupoids, synthesized from the work of many people, 
in particular, Weinstein, Mackenzie, and Ping Xu. 
We then summarize the general theory of multiplicative multivector fields and differential forms on Lie groupoids in Section \ref{multipli}, mainly according to Henrique Bursztyn and Alejandro Cabrera \cite{BC12}. The next two sections deal with the infinitesimal version of multiplicative multivectors and forms on Lie groupoids: they are tensors on Lie algebroids called ``morphic''.
In Section \ref{linearmorphic}, we first review linear vector fields and, more generally, linear multivectors on vector bundles. They define derivations of the graded commutative algebra of sections of their exterior algebra.
We then introduce morphic vector and multivector fields on Lie algebroids. These correspond to derivations of the Gerstenhaber algebra of the Lie algebroid, i.e., the algebra of sections of its exterior algebra equipped with the Schouten--Nijenhuis 
bracket\footnote{The Schouten--Nijenhuis bracket is named after the Dutch mathematician Jan A. Schouten (1883--1971) and his student Albert Nijenhuis (1926--2015). In 1955, four years after he defended his thesis, Nijenhuis vastly extended Schouten's work of 1940 on graded brackets}. We show that the infinitesimal of a multiplicative multivector on a Lie groupoid is a morphic multivector on its Lie algebroid.
In Section~\ref{IM}, we consider the dual case of linear differential forms on vector bundles and morphic differential forms on Lie algebroids, we describe their structure and the correspondence between multiplicative forms on a Lie groupoid and morphic differential forms on its Lie algebroid.
Finally, in Section \ref{GG}, we introduce elements of the {generalized geometry} of manifolds, and we summarize the results on {multiplicative generalized complex structures} on Lie groupoids due to Madeleine Jotz, Mathieu Sti\'enon and Xu~\cite{JSX11}, which unify previous results on symplectic and on holomorphic structures of Lie groupoids.  

\medskip

The importance of Poisson-Lie groups appeared in the theory of integrable systems and their dressing transformations in the 1980s. The successive developments of the theory of multiplicativity have proved that this concept is fundamental.  
Many authors have made important contributions to the subject, and I shall not attempt to cite all of them, but I want
to emphasize the role of our colleagues 
Janusz Grabowski, 
Pawe{\l} Urba{\'n}ski, Katarzyna Grabowska, 
and all their co-authors, 
in Poland and elsewhere, in the development of the theory of Lie algebroids and Lie groupoids and its applications to mechanics, 
and to recall the pioneering role of  {W{\l}odzimierz Tulczyjew} and that of our late colleague and friend, {Stanis{\l}aw Zakrzewski}.

\medskip

This survey is a summary of results where only a few of the proofs are sketched. Complete proofs can be found in the original papers cited below.
Manifolds, tensors and morphisms are assumed to be smooth.
We shall often use the same notation for a morphism between finite-rank vector bundles and its transpose. When introducing local coordinates, we use the Einstein summation convention.

\section{Multiplicative tensors on Lie groups}\label{groups}
We shall deal with the now classical theory of multiplicative multivectors on Lie groups, and we shall see that the need to introduce Lie groupoids arises already at this stage. 
Poisson-Lie groups are named after Sim\'eon-Denis Poisson (1781-1840) and Sophus Lie (1842-1899).
Poisson introduced what became known as the Poisson brackets when, in 1809, he expressed the time derivative of the ``arbitrary constants that became variable'' under a perturbation of the motion of two celestial bodies as a linear combination of the partial derivatives of this perturbation. The Poisson brackets are the coefficients of these expansions and they are time-independent. Lie defined ``function groups'' in 1872 and, in 1888, he derived the 
condition for a contravariant, skew-symmetric tensor of degree 2 to determine a bracket that satisfies the Jacobi 
identity.\footnote{See equation \eqref{Lie} below. For elements of the history of symplectic and Poisson geometry, see \cite{K14}, parts II and III.}  Roughly a century later, Drinfeld imposed a multiplicativity condition on the Poisson structure on a Lie group, defining the ``Hamilton-Lie groups'' \cite{D83} that would later be called ``Poisson-Lie groups'' or ``Poisson groups'', and constitute a semi-classical version of quantum groups.

\subsection{Multiplicative vector fields} 
Consider a {Lie group}, $G$.
By definition, $G$ is a smooth manifold with a group structure such that the multiplication, 
$m: G \times G \to G$, $(g,h) \mapsto m(g,h)= gh$,
is a smooth map, as is the inversion, $\iota:  G \to G$, $g \mapsto \iota(g)= g^{-1}$.

Let $X$ be a vector field on $G$. Can we compare the value of $X$ at $gh$ with the values of $X$ at $g$ and at $h$?
We can left-translate the value of $X$ at $h$ by the tangent of the left translation by $g$,
and 
we can right-translate the value of $X$ at $g$ by the tangent of the right translation by $h$.
Both resulting vectors are tangent to $G$ at $gh$. 
So we can add these two vectors, whence the following definition.
\begin{definition}
A vector field $X$ on a Lie group is  {\it multiplicative} if, for all $g, h \in G$, 
\begin{equation}\label{eq1}
{X_{gh} = g \cdot X_h + X_g \cdot h.}
\end{equation} 
We have abbreviated $(T_h\lambda_g) (X_h)$ to $g \cdot X_h $, where 
$T_h\lambda_g$ is the tangent map at $h\in G$ to $\lambda_g$, the left translation by $g \in G$.
Similarly, 
we have denoted $(T_g\rho_h) (X_g)$ by $ X_g \cdot h$, where $\rho_h$ is the right translation by $h \in G$.
\end{definition}
It is easy to prove the existence of multiplicative vector fields on Lie groups.
Consider the Lie algebra  $\mathfrak g = \mathrm{Lie}(G)$ of $G$. 
For any element $x \in \mathfrak g$,
let $x^\lambda$ be the left-invariant vector field defined by $x$ so that, for all $g\in G$,
$(x^\lambda)_g = g \cdot x$,
and let $x^\rho$ be the right-invariant vector field defined by $x$, so that for all $g\in G$,
$(x^\rho)_g=  x \cdot g$.
\begin{proposition}
For any $x \in \mathfrak g$, the vector field $x^\lambda - x^\rho$ on $G$ is  {{multiplicative}}.
\end{proposition}
\noindent{\it Proof}. The result follows from the properties $\lambda_{gh} = \lambda_{g} \circ \lambda_{h}$
and
$\rho_{gh} = \rho_{h} \circ \rho_{g}$,
and the fact that left- and right-translations commute. \qed

\medskip

A multiplicative vector field that can be written as $X = x^\lambda - x^\rho$ for some $x \in \mathfrak g$ is called \textit{exact} or
a \textit{coboundary}, terms that will be explained in Section \ref{cocycle}.
  
\medskip
 
In the following proposition we describe the flow of a multiplicative vector field, thereby obtaining a justification for the term ``multiplicativity''.
\begin{proposition}
Let $\phi_t$ be the local flow of a vector field $X$ on a Lie group $G$, defined for $t$ in an interval $I$ of $\mathbb R$ containing $0$. 
The vector field $X$ is multiplicative if and only if, for all $g$ and $h \in G$, and all $t \in I$,  
\begin{equation}\label{mult}
{\phi_t(gh)= \phi_t(g) \phi_t(h).}
\end{equation}
\end{proposition}
\noindent{\it Proof}. In fact, if equation \eqref{mult} is satisfied, then $X = \frac{d}{dt}\phi_t|_{t=0}$ satisfies
$$X_{gh} = \frac{d}{dt}\phi_t (g) \phi_t(h)_{|{t=0}} = \frac{d}{dt}\phi_t(g)_{|t=0} \cdot h +  g \cdot \frac{d}{dt}  \phi_t(h)_{|{t=0}},$$ 
and conversely, by integration. \qed

\medskip

In the case of an exact multiplicative vector field, $X = x^\lambda - x^\rho$, the flow of $X$ satisifes 
$\phi_t(g) = \exp(-tx) g \exp(tx)$.

\subsection{Multiplicative multivectors}
The definition of multiplicativity is generalized to contravariant tensor fields of degree $k > 1$ in the obvious way, replacing  
 $T_h\lambda_g$ by the tensor power of degree $k$, $\otimes^k(T_h\lambda_g)$, 
 and replacing $T_g\rho_h$ by $\otimes^k(T_g\rho_h)$. 

A skew-symmetric contravariant tensor field is called a multivector field, or simply a {\it multivector}. If it is of degree $k$, it is called a $k$-vector field, or, simply, a {$k$-vector}.
In particular, a $2$-vector is called a {\it bivector}, and is often 
denoted by $\pi$.
Bivectors play an important role in {Poisson geometry}. 

Any $q \in \wedge^k \mathfrak{g}$ defines an exact multiplicative $k$-vector,
$Q = q^\lambda - q^\rho$.
In particular, an element $r \in \wedge^2 \mathfrak{g}$ defines an exact multiplicative bivector, 
$$\pi = r^\lambda - r^\rho.$$

The notion of ``multiplicativity'' for bivectors and, more generally, for multivectors  appeared in the 1990s, together with the more general notion of {affine multivectors}, in the papers of Jiang-Hua Lu and Weinstein \cite{LW90},  Weinstein \cite{W90}, 
Pierre Dazord and Daniel Sondaz \cite{DS91}, and in \cite{K92}. 

\bigskip

Recall that a smooth manifold, $M$, is a {\it Poisson manifold} if the space of smooth real-valued functions 
over the manifold is a Poisson algebra, i.e., if there is a {skew-symmetric ${\mathbb R}$-bilinear bracket}, $\{~,~\}$, 
on $C^\infty (M)$ which 

(i) is a {derivation} in each argument, and 

(ii) satisfies the {Jacobi identity},
$$
\{f_1,\{f_2,f_3 \}\} + \{f_2,\{f_3,f_1 \}\} + \{f_3,\{f_1,f_2 \}\} = 0,
$$
for all $f_1$, $f_2, f_3 \in C^\infty (M)$. 

Because of condition (i), there is a unique bivector $\pi$ on $M$ such that   
$$
{\pi(df_1,df_2) =\{f_1,f_2\},}
$$ 
for all $f_1$, $f_2 \in C^\infty (M)$. 
Condition (ii)
can then be expressed as the vanishing of a quadratic expression in the components, $\pi^{ij}$, of $\pi$ in a local coordinate system,
\begin{equation}\label{Lie}
{\frac{\partial \pi^{ij}}{\partial x^\ell} \pi^{k\ell} + \frac{\partial \pi^{jk}}{\partial x^\ell} \pi^{i\ell} + \frac{\partial \pi^{ki}}{\partial x^\ell} \pi^{j\ell} = 0.}\quad  \quad 
\end{equation}
In invariant formulation, equation \eqref{Lie} is 
\begin{equation}
[ \pi, \pi ] =0,
\end{equation}
where the bracket is the Schouten--Nijenhuis bracket of multivectors. A bivector such that 
$[ \pi, \pi ] =0$ is called a {\it Poisson bivector}.

\begin{definition} A {\it Poisson-Lie group}, or simply a {\it Poisson group}, is a Lie group equipped with a multiplicative Poisson bivector.
\end{definition} 

The following theorem is an alternative characterization of multiplicative Poisson bivectors.
\begin{theorem}
Let $\pi$ be a Poisson bivector on a Lie group $G$, and let $\{~,~\}$ be the associated Poisson bracket defined by $\pi(df_1,df_2)= \{f_1,f_2\}$, for $f_1, f_2 \in C^\infty (G)$.
Then $\pi$ is multiplicative if and only if the group multiplication, $m: (g,h) \mapsto gh$, is a {Poisson morphism} from $G \times G$ to $G$, where $G \times G$ is equipped with the product Poisson structure.
\end{theorem}
\noindent{\it Proof}. Explicitly, $m : G \times G \to G$ is a Poisson morphism 
if and only if, for all $f_1, f_2 \in C^\infty (G)$ and $g, h \in G$, 
$\{f_1,f_2\}(gh) = \{f_1(g,.),f_2(g,.)\}(h) + \{f_1(.,h),f_2(.,h)\}(g)$,
which is equivalent to the multiplicativity of $\pi$. \qed

\subsection{From QISM to Poisson groups}
The theory of  Poisson groups originated in
the theory of the quantum inverse scattering method (QISM) and dressing transformations, that was the work of the Saint-Petersburg (then Leningrad) school of Ludwig Faddeev, whose best-known members were  P. P. Kulish, Alexei Reiman, Mikhael Semenov-Tian-Shansky, Evgeny Sklyanin, and Leon Takhtajan, 
and in the work of I. M. Gelfand and Irene Ya. Dorfman in Moscow on the relationship between solutions of the classical Yang-Baxter equation\footnote{The various forms of the Yang-Baxter equation are of importance in statistical mechanics, in the theory of integrable systems and in that of quantum groups. The equation is named after Chen Ning Yang, born in 1922, who was awarded the Nobel prize in physics in 1957, was a professor at the State University of New York at Stony Brook, has been emeritus since 1999, and has now returned to China, and
Rodney J. Baxter, born in 1940, whose work is mainly in statistical mechanics, and is professor emeritus at the Australian National University in Canberra.} and Hamiltonian (i.e., Poisson) structures \cite{GD80} \cite{GD82}.

Drinfeld, in his 3-page article \cite{D83}, motivated by Semenov-Tian-Shansky's ``What is a classical r-matrix?'' \cite{S83}, 
introduced Lie groups with a Hamiltonian (i.e., Poisson) structure satisfying the requirement that the group multiplication be a Poisson map
from the product Poisson manifold $G \times G$ to $G$, which he called {\it Hamilton-Lie groups}.\footnote{In a lecture in Oberwolfach in the summer of 1986, published the following year \cite{K87}, I called Drinfeld's Hamilton-Lie groups ``Poisson-Drinfeld groups'', but the name was not retained by later authors. }
The Hamilton-Lie groups of Drinfeld gained importance when he showed that they appeared as the semi-classical version of the quantum groups which he would introduce three years later in his address \cite{D86} to the ICM in Berkeley in 1986. There, he defined them as those groups equipped with a Poisson bracket which is ``compatible with the group operation'', and he called them {\it Poisson-Lie groups}.
This was the term adopted by Lu and Weinstein in their influential article \cite{LW90} and by Izu Vaisman \cite{V94}. 
Later, it was often shortened to {\it Poisson groups}, which is how Semenov-Tian-Shansky had named 
them \cite{S85}. We refer to them as Poisson groups in the following pages.
\subsection{Multiplicativity as a cocycle condition}\label{cocycle}
The cohomological approach, already present in Drinfeld's paper and in \cite{LW90}, was developped in \cite{K92}. 
For any field, $Q$, of multivectors of degree $k$ on a Lie group, $G$, define the map  
$\rho(Q): G \to \wedge^k \mathfrak{g}$ by
$$\quad \rho(Q)(g)= Q_g \cdot g^{-1},$$ 
and the map $\lambda(Q): G \to \wedge^k \mathfrak{g}$ by
$$\quad \lambda(Q)(g)= - g \cdot Q_{g^{-1}}.$$
\begin{theorem}\label{group_cocycle}
The following properties are equivalent.

\noindent (i) $Q$ is multiplicative,

\noindent (ii) $\rho(Q)$ is a $1$-cocycle of $G$ with values in $\wedge^k \mathfrak{g}$,

\noindent (iii) $\lambda(Q)$ is a $1$-cocycle of $G$ with values in $\wedge^k \mathfrak{g}$.
\end{theorem}
\noindent{\it Proof}. If $Q$ is multiplicative, then 
$$\rho(Q)(gh)= Q_g \cdot h \cdot (gh)^{-1} + g \cdot Q_h \cdot (gh)^{-1} = \rho(Q)(g)+ g \cdot \rho(Q)(h) \cdot g^{-1}.$$ 
Thus $\rho(Q)$ is a {$1$-cocycle of the group $G$ with values in $\wedge^k \mathfrak{g}$,
where $G$ acts by the adjoint action. 
The converse follows from the same calculation. 
The computation for $\lambda(Q)$ is similar. \qed
  
\medskip

If a multivector $Q$ is multiplicative and exact, then $\rho(Q)= \lambda (Q)$. 
In fact, if $Q = q^\lambda - q^\rho$, then $\rho(Q)(g) = g \cdot q \cdot g^{-1} - q = \lambda (Q)(g)$.
A $k$-vector is multiplicative and exact if and only if $\rho(Q)$ (or $\lambda(Q)$) is a {1-coboundary} of $G$ with values in $\wedge^k \mathfrak{g}$.
These facts justify the terms {\it exact} and {\it coboundary}.

When a bivector $\pi$ on $G$ is equal 
to $r^\lambda - r^\rho$, where $r \in \wedge^2{\mathfrak{g}}$, it is multiplicative. It remains to express the fact that it is a Poisson bivector.
Let $[~,~]$ be the {\it algebraic Schouten bracket} defined as the extension of the Lie bracket of $\mathfrak{g}$ as a biderivation of
the exterior algebra, $\wedge^\bullet {\mathfrak{g}}$. (It is also the trace on $\mathfrak{g}$ of the Schouten--Nijenhuis bracket of left-invariant multivectors on the Lie group, $G$.) For $r \in \wedge^2{\mathfrak{g}}$, $[r,r]$ is an element in $\wedge^3{\mathfrak{g}}$.

\begin{theorem}
For $\pi = r^\lambda - r^\rho$ to be a {Poisson bivector} it is necessary and sufficient that the element $[r,r] \in \wedge^3{\mathfrak{g}}$ be $Ad_G$-invariant. A sufficient condition is  $[r,r] = 0$. 
\end{theorem}

If $[r,r]$ is $ad_{\g}$-invariant, $r$ is said to be a solution of the {\it generalized classical Yang-Baxter equation}. The condition $[r, r] = 0$ is called the {\it classical Yang-Baxter equation}.
A~Poisson group, defined by $r \in \wedge^2{\mathfrak{g}}$ satisfying the classical Yang-Baxter equation, is said to be {\it triangular}.

\subsection{An equivalent definition of multiplicativity for vector fields}
We first recall the definition of the {\it tangent group} of a Lie group.
The tangent bundle $TG$ of a Lie group $G$ is a Lie group where the multiplication is the tangent of the multiplication in $G$ and the inversion is the tangent of the inversion in $G$. In other words the multiplication, $\times_{(TG)} : TG \times TG \to TG $, is such that for $X \in T_gG$ and $Y \in T_hG$, 
\begin{equation}\label{tangent}
X \times_{(TG)} Y =  g \cdot Y +  X \cdot h
\end{equation}
(recall that $g \cdot Y +  X \cdot h  = (T_h\lambda_g)(Y_h) + (T_g\rho_h)(X_g)$). 
The inverse of $X \in T_gG$ is $X^{-1} = (T_g \iota)(X) = - g^{-1} \cdot  X \cdot g^{-1} \in T_{g^{-1}}G.$
\begin{theorem}
A vector field on $G$, $X: G \to TG$, is multiplicative if and only if $X$ is a morphism of groups. 
\end{theorem}
\noindent{\it Proof}. Vector field $X$, seen as a map from $G$ to $TG$, is a morphism of groups if and only if
$X_{gh}= X_g \times_{(TG)} X_h.$
By the definition of multiplication in the group $TG$,
this condition coincides with condition (\ref{eq1}), which is the condition for $X$ to be multiplicative. \qed

\subsection{Questions}
One may ask whether the cotangent bundle of a Lie group can be made into a group. For an interesting answer, we need to introduce groupoids.
Among the main references for the early theory of groupoids, and Lie groupoids in particular, are the works of Ronald Brown, Charles Ehresmann, Jean Pradines, and Paulette Libermann (see the references in \cite{M87} \cite{M05}).

Several questions then arise. Are there multiplicative multivectors on Lie groupoids?
Are there multiplicative differential forms on Lie groups? on Lie groupoids?

We remark immediately that any multiplicative multivector $Q$ on a Lie group vanishes at the identity $e$ of the group, because setting $g = h = e$ in the defining equation 
yields $Q_e = 0$. Since a multiplicative bivector on a Lie group cannot be everywhere non-degenerate, there is no corresponding $2$-form, which implies that there are no multiplicative,
non-degenerate $2$-forms (``multiplicative symplectic forms'') on Lie groups. 

\section{Multiplicativity and groupoids}\label{groupoids}
We shall introduce the concept of groupoids and define the groupoid structure of their tangent and cotangent bundles, and then
explain the role of groupoids in the theory of Poisson groups. 
\subsection{Groupoids generalize groups}
In a {Lie groupoid}, $\varGamma$ and $ M$ are smooth manifolds. There is a {\it source} $\alpha :\varGamma {\to} M$ and a {\it target} $\beta :\varGamma {\to} M$  
which are surjective submersions. 
$$ \begin{array}{rcl}  & \varGamma &   \\  \alpha & \downarrow 
 \downarrow & \beta \\ & M &
 \end{array}$$
There is a partially defined associative {\it multiplication}: 
if $g \in \varGamma$, $h \in \varGamma$, 
then $gh$ is defined if and only if $\alpha(g) = \beta(h) $. 
Let $\varGamma^{(2)}$ be the submanifold of composable pairs in $\varGamma \times  \varGamma$.
The multiplication map $(g,h) \in \varGamma^{(2)} \mapsto gh \in \varGamma$ is smooth.
Each element has an {\it inverse} and the inversion is a smooth map, $\iota : \varGamma \to \varGamma$.

The space $M$, called the {\it base} of the groupoid, is identified with the submanifold of the {unit} elements of the multiplication.

One can think of the elements of $\varGamma$ as paths whose multiplication is defined if the target of one path coincides with the source of the second path.

\begin{example}
The pair groupoid on base $M$ is $M \times M$ with the projections, 
$\alpha=pr_2$, $\beta = pr_1$, and multiplication map, {$(x,y)(y,z) = (x,z)$.}
The inversion is $(x,y) \mapsto (y,x)$.
The base $M$ is identified with the diagonal submanifold, $\{(m,m) \in  M \times M | m \in M \}$.
\end{example}
\begin{example}
 Any Lie group can be considered as a groupoid over a point.
\end{example}
\begin{example}\label{gaugegpd}
The gauge groupoid of a principal bundle $P \to M$, whose structure group is a Lie group $G$, is a Lie groupoid defined as follows\footnote{The Lie groupoid associated to a 
principal bundle was defined by Ehresmann in his lecture at the Colloque de topologie in Brussels in 1950. See Libermann [L 69] [L71].}.  
Let $\varGamma$ be the quotient of $P \times P$ by the diagonal action of $G$ on the right: 
$P \times P \times G \to P \times P$,
$((u,v), g) \mapsto (ug,vg)$, and let $cl(u,v)$ denote the equivalence class of $(u,v) \in P \times P$. The source (resp., target) is induced by the composition of the second (resp., first) projection from $P \times P$ to $P$ composed with the projection $P \to M$. The product of two composable elements in the quotient $\varGamma$ is defined by choosing representatives
$(u,v)$ and $(v,w)$ and setting                                                   
$cl(u,v) \, cl(v,w) = cl(u,w)$, thus the inverse of $cl(u,v)$ is $cl(v,u)$. Then $\varGamma$ is a Lie groupoid with base $M$, the quotient\footnote{For the notion of the quotient of a groupoid by a fibration, see [M05].} of the pair groupoid 
$P \times P $ with base $P$ by the fibration $P \times P \to \varGamma$ over the projection $P \to M$.
\end{example}

\subsection{The cotangent groupoid of a Lie group}
The tangent bundle of a Lie group is a Lie group, and, like all Lie groups, it can be considered to be a Lie groupoid over a point. What about the cotangent bundle of a Lie group?

For any Lie group, $G$, the cotangent bundle $ T^*G$ has a canonical structure of a Lie groupoid with base ${\mathfrak g}^*$,
$$ \begin{array}{ccc}  & T^*G &   \\  \alpha_{T^*G} & \downarrow 
 \downarrow & \beta_{T^*G} \\ & {\mathfrak g}^* &
 \end{array}$$
For $\xi \in T_g^*G$, let $\alpha_{T^*G} (\xi)=  \xi \circ T_e\lambda _g$ and $\beta_{T^*G} (\xi)=  \xi \circ T_e\rho _g$.
If $\xi  \in T_g^*G$ and $\eta \in T_h^*G$, their product is defined if and only if 
$$\xi \circ T_e\lambda_{g}= \eta \circ T_e\rho_{h},$$
and then it is 
$$\xi \times_{( T^*G)} \eta =  \xi \circ T_e\rho_{h^{-1}}= \eta \circ T_e\lambda_{g^{-1}} \in T_{gh}^*G.$$
This multiplication satisfies, for all $X  \in T_g G$ and $Y \in T_h G$,
\begin{equation}\label{AAA}
(\xi \times_{( T^*G)} \eta) (X \times_{( TG)} Y) = \xi(X) + \eta(Y).
\end{equation}
The inverse of $\xi \in T^*_g G$ is $\xi^{-1}= - \, {}^t(T_{g}\iota) (\xi) = - \, {}^t(T_{e}\lambda_{g^{-1}}) \circ {}^t(T_{g}\rho_{g^{-1}}) (\xi) \in T^*_{g^{-1}}G$.
\subsection{Multiplicative bivectors as morphisms of groupoids} While Poisson groups are characterized by the condition that the Poisson bivector be multiplicative,
they can also be characterized by a morphism property.
If $\pi$ is a bivector on a manifold, $P$, the map  we shall denote by $\pi^\sharp : T^*P \to TP$ is defined by $\pi^\sharp(\xi)= i_\xi \pi = \pi(\xi, \cdot)$, for all $\xi \in T^*P$. We shall often simply write $\pi$ for the map $\pi^\sharp$. Similarly, if $\omega$ is a $2$-form on a manifold, $P$, we shall denote by $\omega^\flat$, or, simply, $\omega$, the map $TP \to T^*P$ defined by $\omega^\flat (X) = - i_X \omega = \omega (\cdot, X)$, for all $X \in TP$.  Then, if $\omega(\pi^\sharp \xi, \pi^\sharp \eta) = \pi(\xi,\eta)$, for all $\xi, \eta \in T^*P$, and if $\pi^\sharp$ is invertible, $\omega^\flat$ is its inverse.
\begin{theorem}\label{2.1}
A Poisson bivector $\pi$ on a Lie group $G$ is multiplicative if and only if $\pi^\sharp$ is a {morphism from the Lie groupoid $T^*G$ over ${\mathfrak g}^*$
 to the Lie group $TG$ considered as a Lie groupoid over a point},
$$ \begin{array}{ccccc}   T^*G & \stackrel{\pi^\sharp}{\to} & TG   \\    \downarrow 
 \downarrow & &
 \downarrow  \\  {\mathfrak g}^* & \to & \{pt\}
 \end{array}.$$
\end{theorem}
This theorem \cite{M92} -- a re-formulation of the multiplicativity property of Poisson bivectors in terms of morphisms of groupoids --
permitted vast generalizations of the concept of multiplicativity.

\section{Poisson groupoids}\label{PG}

In order to generalize the definition of Poisson-Lie groups to 
groupoids, the most straightforward method is to consider a groupoid analogue of the characterization of multiplicativity for Poisson bivectors on Lie groups given in Theorem \ref{2.1}, and then to use it as the definition of a ``multiplicative'' 
Poisson structure on a Lie groupoid. We first need to recall the definition of the Lie algebroid of a Lie groupoid due to Pradines \cite{P67} (see, e.g., \cite{M08}) and to 
describe the Lie groupoid structure of the tangent 
and cotangent bundles of a Lie groupoid.

\subsection{The Lie algebroid of a Lie groupoid}\label{FFFF}
\begin{definition} A real vector bundle, $A \to M$, for which the module of sections has a Lie-algebra structure over the reals and there is a vector-bundle morphism called the {\it anchor}, $\rho: A \to TM$, that satisfy the Leibniz rule,
$$
[X, fY] = f [X,Y] + \rho(X)f \, Y,
$$
for all sections $X$ and $Y$ of $A \to M$ and all $f \in C^\infty (M)$, is called a 
{\it Lie algebroid}.
\end{definition}
\begin{example} Any Lie algebra is a Lie algebroid whose base is a point. In this case the anchor vanishes.
\end{example} 
\begin{example} The tangent bundle $TM$ of any smooth manifold $M$ is a Lie algebroid with base $M$ whose anchor is the identity map of $TM$.
\end{example} 
\begin{example} The cotangent bundle of a Poisson manifold, $(M,\pi)$, is, in a canonical way, a Lie algebroid, where the Lie bracket of exact $1$-forms, $\d f$ and $\d g$, is the de Rham differential of the Poisson bracket of $f$ and $g$. An explicit form of this bracket has been known since 1982 \cite{F82} \cite{MM84} \cite{CDW87}  \cite{H90}. The anchor of the Lie algebroid $T^*M$ is $\pi^\sharp : T^*M \to TM$.
\end{example}  
\medskip

Just as each Lie group has a Lie algebra, each Lie groupoid, $\varGamma$, has a {Lie algebroid}, denoted by $A(\varGamma)$ or, simply, $A$. It is the {vector bundle}, over the base $M$ of $\varGamma$, whose fibers are the tangent spaces to the $\alpha$-fibers of $\varGamma$ at all points of $M$.
Given a section $X$ of $A$, it can be extended uniquely into a right-invariant vector field on $\varGamma$.
By definition, the {\it Lie bracket} of sections of $A$ is the restriction to $M$ of the Lie bracket of the corresponding right-invariant vector fields.
The map $\rho : A \to TM$ that associates to any element of $A$ its projection on $TM$ along the $\beta$-fibers is the 
{anchor} of the Lie algebroid.  
\begin{example} The Lie algebroid of a pair groupoid $M \times M$ is the tangent bundle of $M$.
\end{example} 
\begin{example}\label{gaugealgd}
 The Lie algebroid of the gauge groupoid of a principal bundle $P \to M$ with structure group $G$ is called the ``gauge algebroid'' or the ``Atiyah algebroid'' of the principal bundle\footnote{The exact sequence of Lie algebroids attached to a principal bundle was 
defined by Michael Atiyah in his study of connections on complex 
analytic vector bundles [A57], but he did not consider the associated 
groupoids. Libermann \cite{L69} \cite{L71} identified the Lie algebroid 
of the gauge groupoid of a principal bundle with the ``Atiyah 
algebroid'', calling its elements ``infinitesimal displacements'', a 
notion that had been introduced by Ehresmann in his lecture at the Colloque de topologie in Brussels in 1950.}. Its underlying vector bundle is the quotient of $TP$ by the action of $G$ defined as the tangent prolongation of the action of $G$ on $P$ on the right. The bracket of sections is induced from the Lie bracket of right-invariant vector fields on $TP$.
The anchor is induced from $TP \to TM$, the tangent of the projection $P \to M$.
\end{example} 
\subsection{The dual of a Lie algebroid: linear Poisson structures}

When $A \to M$ is a Lie algebroid, the dual vector bundle, $A^*\to M$, has a {\it linear Poisson structure}.
In fact, there is a Poisson bracket on the space of functions on $A^*$ characterized by the following properties. 

(i) the bracket of functions linear on the fibers of $A^*$ is a function linear on the fibers, defined as the Lie-algebroid bracket of the corresponding sections of $A$, 

(ii) the bracket of functions pulled back to $A^*$ from functions on the base vanishes, and 

(iii) the bracket of a function linear on the fibers of $A^*$, defined by a section $X$ of $A$, and a (pulled back) function on the base is the (pull-back of the) Lie derivative of that function by $\rho(X)$, where $\rho$ is the anchor of $A$.  

\medskip

Conversely, the dual of a vector bundle with a linear Poisson structure is a Lie algebroid \cite{W88}.
 
\subsection{The tangent groupoid of a Lie groupoid}\label{tangentG} 
The construction of the tangent group of a Lie group can be extended to the construction of the tangent groupoid of a Lie groupoid, $\varGamma$,
$ \begin{array}{c} 
T\varGamma    \\  \downarrow 
 \downarrow \\  TM 
 \end{array}$.
The groupoid multiplication in the tangent bundle, $T\varGamma$, of $\varGamma$ is the tangent of the multiplication in $\varGamma$, the inversion is the tangent of the inversion of $\varGamma$, and the source and target are also obtained by applying the tangent functor to the source and target of $\varGamma$. 
 
\subsection{The cotangent groupoid of a Lie groupoid}\label{cotangentG}
For any Lie groupoid, $\varGamma$, the cotangent bundle, $T^*\varGamma$, has the structure of a Lie groupoid with base $A^*$, the dual of the Lie algebroid $A$ of $\varGamma$, $ \begin{array}{ccc}  & T^*\varGamma &   \\ \alpha_{T^*\varGamma}  & \downarrow 
 \downarrow & \beta_{T^*\varGamma} \\ & A^* &
 \end{array}$, defined as follows \cite{CDW87} \cite{M13}.  
 
Let $g$ be an element of $\varGamma$.
The {source} of an element $\xi \in T_g^*\G$ is the element, $\alpha_{T^*\varGamma} \xi$, of $A^*_{\alpha g}$  defined by 
$$\alpha_{T^*\varGamma} \xi(X) = \xi\big ( g  \cdot  ( X - \rho (X)) \big ),$$ for all $X\in A_{\alpha g}$ (where $\rho : A \to TM$ is the anchor of $A$, and $T_{\alpha g}\lambda_g$ was denoted by $g \, \cdot$),

The {target} of $\xi$ is the element, $\beta_{T^*\varGamma} \xi$, of $A^*_{\beta g}$ defined by 
$$\beta_{T^*\varGamma} \xi(X) = \xi(X \cdot g)$$ for $X\in A_{\beta g}$ (where $(T_{\beta g}\rho_g)(X)$ was denoted by $X \cdot g$).

Let $g, h \in \varGamma$.
If $\xi \in T_g^*\G$ and $\eta \in T_h^*\G$, their {product} is defined if and only if $\alpha_{T^*\varGamma}\xi = \beta_{T^*\varGamma} \eta$, and then $\alpha g = \beta h$ and the product is the element $ \xi \times_{(T^*\varGamma)} \eta \in T_{gh}^*\G$,
such that
\begin{equation}\label{cotangent}
(\xi \times_{(T^*\varGamma)} \eta) (X \times_{(T\varGamma)} Y)= \xi (X) +  \eta(Y),
\end{equation} 
for all $X \in T_{g}\varGamma$ and  $Y\in T_{h}\varGamma$. 

The inverse of $\xi \in T^*_g \varGamma$ is the element in $T^*_{g{-1}} \varGamma$ opposite of the image of $\xi$ under the transpose of $T_{g}\iota$.
  
\medskip
  
The groupoid structure of $T^*\varGamma$ can be obtained as a particular case of Pradines's general theory of duality for ``vector-bundle groupoids'' \cite{P88}. 

\subsection{Definition of Poisson groupoids}\label{3.4}

\medskip

\begin{definition} 
A {\it Poisson groupoid} is a Lie groupoid equipped with a Poisson bivector that defines a map, $\pi^\sharp : T^*\varGamma \to T\varGamma$, over a map, $  {\underline \pi} : A^* \to TM$, which is a {morphism of Lie groupoids},
$$ \begin{array}{ccc}  T^*\varGamma &  \stackrel{\pi^\sharp}{\to} & T\varGamma   \\    \downarrow  \downarrow &  &   \downarrow 
 \downarrow \\ {A}^* &  \stackrel{\underline \pi}{\to} &   TM
 \end{array}.$$
\end{definition}

A submanifold of a Poisson manifold is called {\it coisotropic} if the Poisson bracket of two functions that vanish on the submanifold also vanishes on this submanifold.
Equivalently, a submanifold of a Poisson manifold is {coisotropic} if and only if, at each point in the submanifold, the image by the Poisson map of the orthogonal of the tangent space is contained in the tangent space. 
The original definition of Poisson groupoids was given in 1988 by Weinstein in \cite{W88} as a generalization of both the {Poisson-Lie groups} and the {symplectic groupoids}.
They were defined as {Lie groupoids with a Poisson structure}, $\pi$, which
is {\it multiplicative} in the sense that the graph of the partially defined multiplication, $m: \varGamma^{(2)} \subset \varGamma \times \varGamma \to \varGamma$, is a {coisotropic} submanifold of the product Poisson manifold, 
$\varGamma_\pi \times \varGamma_\pi \times \varGamma_{-\pi}$. 
Claude Albert and Dazord \cite{AD90}, then Mackenzie and Xu \cite{MX94} proved that the {definition of Poisson groupoids in terms of morphisms of groupoids is equivalent} to Weinstein's original definition.

To see that Poisson groupoids, defined by the coisotropy of the graph of their multiplication, generalize Poisson groups, recall that a Poisson bivector on a Lie group, $G$, is multiplicative if and only if the multiplication is a Poisson map from $G \times G$, with the product Poisson structure, to $G$ and apply the following lemma \cite{W88}.
\begin{lemma} Let $\pi$ be a Poisson bivector on a Lie group $G$.  A map from  $G_\pi  \times G_\pi $ to $G_\pi $ is Poisson if and only if its graph is coisotropic in 
$G_\pi \times G_{\pi} \times G_{-\pi}$.
\end{lemma}

\subsection{Symplectic groupoids as Poisson groupoids}
Recall that, by definition, a symplectic structure on a manifold is a non-degenerate, closed $2$-form. (The manifold is then necessarily of even dimension.)
A manifold is symplectic if and only if it is a Poisson manifold with a non-degenerate Poisson bivector.

A submanifold of a symplectic manifold is called {\it Lagrangian} if it is coisotropic and of minimal dimension (half the dimension of the ambient symplectic manifold).
Equivalently, a submanifold of a symplectic manifold is {Lagrangian} if and only if the symplectic orthogonal of each tangent space is equal to that tangent space.

{\it Symplectic groupoids} were originally defined as Lie groupoids with a symplectic structure, $\omega$, which is {\it multiplicative}, in the sense that the graph of the multiplication, $m$, is a 
{Lagrangian} submanifold of the symplectic 
manifold $\varGamma_\omega \times \varGamma_{\omega} \times \varGamma_{-\omega}$ \cite{K86} \cite{KM93} \cite{W87} \cite{Z90}. 
It follows from the definitions that {symplectic groupoids coincide with those Poisson groupoids whose Poisson structure is non-degenerate}.
Although it is particularly important and has many features not present in the Poisson case in general, the case of symplectic groupoids can thus be treated as a particular case of the more general theory of Poisson groupoids.

It is easy to see that a 2-form, $\omega$, on a Lie groupoid is multiplicative if and only if it satisfies the condition,
$$
m^*(\omega) = pr_1^*(\omega) +  pr_2^*(\omega),
$$
where $pr_1$ (resp., $pr_2$) is project$  $ion to $\varGamma$ of the first (resp., second) factor in the domain of the multiplication, $m: \G^{(2)} \subset \G \times \G \to \G$. 
It was then 
proved that $\omega$ is multiplicative if and only if $\omega^\flat: T\varGamma \to T^*\varGamma$
is a morphism of Lie groupoids, where $T\varGamma$ and $T^*\varGamma$ are equipped with the groupoid structures over 
$TM$ and $A^*$, respectively, defined in Sections \ref{tangentG} and \ref{cotangentG}.

\begin{remark} 
As we remarked above, there are 
{no multiplicative symplectic forms on Lie groups}.
On the other hand, in the 1980s, Andr\'e Lichnerowicz, together with Alberto Medina, studied Lie groups equipped with a 
{left-invariant symplectic form} and called them ``symplectic groups'' (see Lichnerowicz's lecture at the Colloque Souriau of 1990 \cite{L91}). Despite the terminology, these ``symplectic groups''  are not a particular case of the symplectic groupoids defined above.
\end{remark} 

\subsection{The infinitesimal of a Poisson groupoid}
The infinitesimal of a Poisson group is a {Lie bialgebra}. Mackenzie and Xu \cite{MX94} introduced the definition of Lie bialgebroids and showed that the infinitesimal of a Poisson groupoid is a {\it Lie bialgebroid}, i.e., a pair of Lie algebroids in duality, $(A,A^*)$, such that the differential on the sections of the exterior algebra of $A^*$, defined by the Lie-algebroid structure of $A$, is a derivation of the Schouten--Nijenhuis bracket that extends the Lie-algebroid bracket of $A^*$ (see \cite{K95}).

In a Poisson groupoid, $\varGamma$, the base is a coisotropic submanifold. The dual, $A^*$, of the Lie algebroid, $A$, of $\varGamma$ inherits a 
Lie-algebroid structure because it is isomorphic to the conormal bundle of a coisotropic submanifold\footnote{For the Lie-algebroid structure on the conormal bundle of a coisotropic submanifold of a Poisson manifold, see \cite{W88}.} of the Poisson manifold $\varGamma$.
From the multiplicativity of the Poisson structure of $\varGamma$, it follows that the pair $(A,A^*)$ is a {Lie bialgebroid} \cite{MX94}. Other proofs of this fact were given subsequently in \cite{MX98} and by David Iglesias-Ponte, Camille Laurent-Gengoux and Xu in \cite{ILX12} (see also \cite{LSX11}).

\subsection{The infinitesimal of a symplectic groupoid}

In a symplectic groupoid, $\varGamma \overrightarrow{\rightarrow} M$, the base is a Lagrangian submanifold. The symplectic form of $\varGamma$ induces a unique Poisson structure on the base such that the source map is a Poisson map and the 
target map is anti-Poisson, i.e., it reverses the sign of the Poisson bracket. 
The Lie algebroid, $A$, of $\varGamma$ is isomorphic to the cotangent bundle of the base, $T^*M$, equipped with the Lie-algebroid structure coming from the induced Poisson structure on 
$M$ \cite{CDW87},\footnote{In fact, it is by transporting the Lie-algebroid bracket on the space of sections of $A$ to a bracket on the space of $1$-forms on $M$ by means of this isomorphism of vector bundles that the authors defined the bracket of $1$-forms on a Poisson manifold, independently of Benno Fuchssteiner \cite{F82} and of Franco Magri and Carlo Morosi \cite{MM84} who had introduced it earlier.}
and the pair $(A,A^*)$ is isomorphic to the Lie bialgebroid $(T^*M,TM)$. It is the non-degeneracy assumption on the Poisson structure of $\varGamma$ that implies that $A$ and $T^*M$ are isomorphic as vector bundles. It is still true in the more general case of a Poisson groupoid that the base manifold inherits a Poisson structure but, in  general, the Lie bialgebroids $(A,A^*)$ and $(T^*M,TM)$ are distinct. 
To summarize, in a Poisson groupoid, the base is coisotropic and a Poisson manifold, and $A^*$ is isomorphic to the conormal bundle $\nu ^* M$ of the base, which is a Lie algebroid because the base is coisotropic. In a symplectic groupoid, the base is Lagrangian and a Poisson manifold, and $A \simeq \nu M \simeq T^*M$ and $A^* \simeq  \nu^* M \simeq TM$.

Conversely, conditions for the existence of a symplectic groupoid ``integrating'' a given Poisson manifold are to be found in an important paper by Marius Crainic and Rui Fernandes published in 2004 \cite{CF04}.



\section{Multiplicative forms and multivectors on Lie groupoids}\label{multipli}

We have defined multiplicative multivectors on Lie groups,
multiplicative bivectors on Lie groupoids and multiplicative $2$-forms on Lie groupoids.
Here, we shall define, more generally, {multiplicative $k$-forms} and {$k$-vectors} on Lie groupoids, for $k$ a nonnegative integer.
 
\subsection{Multiplicative functions, 1-forms and vector fields on Lie groupoids}\label{4.1}
Multiplicative $1$-forms and vector fields on Lie groupoids were discussed in \cite{MX98}.

A real-valued function, $f$, on a Lie groupoid, $\varGamma$, is said to be {\it multiplicative} if it defines a morphism of groupoids from 
$\varGamma$ to $\mathbb R$, which is here considered to be a Lie groupoid over a point, i.e.,
for $g$ and $h$ composable elements of $\varGamma$, 
$$
f(gh) = f(g) + f(h).
$$

A 1-form, $\omega$, on a Lie groupoid, $\varGamma$, is said to be {\it multiplicative} if it satisfies 
$${m^*(\omega) = pr_1^*(\omega) +  pr_2^*(\omega)},$$
 or, equivalently, if it is a morphism of groupoids from $\varGamma$ as a groupoid over $M$, to $T^*\varGamma$ as a groupoid over $A^*$, i.e.,
for all $g$ and $h$ composable in $\varGamma$, 
$$
\omega_{gh} = \omega_g \times_{(T^*\varGamma)} \omega_h, 
$$ 
or, equivalently, if it is a morphism of groupoids from $T\varGamma$ as a groupoid over $TM$, to $\mathbb R$ considered as a groupoid over a point, i.e.,
for all $X$ and $Y$ composable in $T\varGamma$, 
\begin{equation}\label{1form}
\omega(X \times_{(T\varGamma)} Y)= \omega(X) + \omega(Y).
\end{equation} 

A multiplicative vector field is defined as a groupoid morphism from $\varGamma$ as a groupoid over $M$ to  $T\varGamma$ as a groupoid over $TM$.

The Lie bracket of two multiplicative vector fields is multiplicative.
The contraction of a multiplicative 1-form with a multiplicative vector field is a multiplicative function.

\subsection{Multiplicative 2-forms on Lie groupoids}\label{4.2} We have seen above that multiplicative 2-forms on a Lie groupoid $\varGamma \overrightarrow{\rightarrow} M$ appeared when Weinstein defined symplectic groupoids in 1987 \cite{W87} (but he did not use the term ``multiplicative'' in this context). 
Here we summarize several characterizations of {\it multiplicative $2$-forms} (see \cite{MX00}). 
\begin{theorem}\label{2forms}
Let $\omega$ be a $2$-form on a Lie groupoid $\varGamma$ over $M$, with multiplication $m$, and let $A$ be the Lie algebroid of $\Gamma$. The following properties are equivalent.

\noindent(i) The graph of $m$ is Lagrangian in $\varGamma_\omega \times \varGamma_{\omega} \times \varGamma_{-\omega}$.

\noindent(ii) ${m^*(\omega) = pr_1^*(\omega) +  pr_2^*(\omega)}$.

\noindent(iii) $\omega^\flat : T \varGamma \to T^* \varGamma$ is a {morphism of Lie groupoids}, over a vector-bundle morphism ${\underline \omega} : TM \to A^*$,
$$
\begin{array}{ccc}  T \varGamma &  \stackrel{\omega^\flat}{\to} & T^* \varGamma   \\    \downarrow  \downarrow &  &   \downarrow 
 \downarrow \\
TM & \stackrel{\underline \omega} {\to} &  A^*
\end{array},
$$
i.e., for any pair $(X,Y)$ of composable elements in $T\varGamma$,
\begin{equation}\label{forms}
\omega^\flat(X \times_{(T\varGamma)} Y) = \omega^\flat(X) \times_{(T^*\varGamma)} \omega^\flat(Y).
\end{equation}
\noindent(iv) The 2-form $\omega$ defines a morphism from the direct sum Lie groupoid $T\varGamma \oplus T\varGamma$ as a groupoid over $TM \oplus TM$, to the additive group
${\mathbb R}$ considered as a groupoid over a point.\\ 
\noindent(v) The graph of $\omega^\flat :  T \varGamma \to T^* \varGamma$ is a {Lie subgroupoid} of the direct sum Lie groupoid, $T\varGamma \oplus T^*\varGamma \overrightarrow{\rightarrow}  TM \oplus A^*$, over some vector sub-bundle of $TM \oplus A^*$.
\end{theorem}
We shall prove that (iii) implies (iv). 
Applying both sides of formula (\ref{forms}) for $X=X_1$ and $Y = Y_1$ to the product of two composable elements, $X_2$ and $Y_2$ in $T\varGamma$, and using equation (\ref{cotangent}), we obtain 
$$  \omega^ \flat (X_1\times_{(T\varGamma)} Y_1)(X_2 \times_{(T\varGamma)} Y_2) = (\omega ^ \flat (X_1) \times_{(T^*\varGamma)} \omega^ \flat (Y_1))(X_2 \times_{(T\varGamma)} Y_2) 
$$
$$
=\omega^ \flat (X_1)(X_2) + \omega^ \flat (Y_1)(Y_2)
$$
and therefore 
\begin{equation}\label{sum2forms}
\omega(X_1 \times_{(T\varGamma)} Y_1, X_2 \times_{(T\varGamma)} Y_2) =  \omega(X_1,X_2) + \omega(Y_1, Y_2),
\end{equation}
which proves (iv). The converse follows by the same computation.

\medskip
In particular, a {\it presymplectic groupoid}~\cite{BCWZ04} is a Lie groupoid, $\G$, equipped with a  multiplicative, closed 2-form $\omega$  satisfying a non-degeneracy assumption, that is milder than being non-degenerate, which is the condition that the map from the Lie algebroid, $A$, of $\varGamma$ to $TM \oplus T^*M$, $X  \mapsto (T\beta(X), \omega^\flat (X)|_{TM})$ be injective. The $2$-form $\omega$ is identified with the map, $\omega^\flat : T\G \to T^*\G$, over a vector-bundle morphism, ${\underline \omega} : TM  \to A^*$.
The transpose of ${\underline \omega}$ is a vector-bundle morphism from $A$ to $T^*M$, denoted by $\rho^*_\omega$ 
in~\cite{BCWZ04}, then by $\sigma$ in \cite{BCO09} and called there an IM-form. See Section \ref{IMform} below.

In Section \ref{ortiz} we shall give a further characterization of multiplicative, closed 
2-forms on $\varGamma$ as Dirac structures in the standard Courant algebroid, $T\G \oplus T^*\G$.
\subsection{Multiplicative $k$-forms on Lie groupoids}
It is straightforward to extend the definition of multiplicative $1$-forms and $2$-forms to the case of forms of higher degree \cite{BCO09} \cite{ BC12}. Let $k$ be a positive integer. By definition, on a Lie groupoid, $\varGamma$, with multiplication $m$, a $k$-form $\omega$  is {\it multiplicative} if 
$$
m^*(\omega) = pr_1^*(\omega) +  pr_2^*(\omega).
$$
Multiplicative $k$-forms, for $k \geq 1$, can be characterized as morphisms of groupoids as follows.
\begin{theorem} 
A $k$-form $\omega$ on a Lie groupoid, $\G$, with Lie algebroid $A$ is multiplicative if and only if it defines a morphism of Lie groupoids, 
$$ \begin{array}{ccc} \oplus^{k-1} T \G &\stackrel{}{\to} & T^*\G \\ \downarrow \downarrow &
 &\downarrow  \downarrow \\  \oplus^{k-1} TM & \stackrel{}{\to} & A^* \\
 \end{array},$$
or, equivalently,
 $$ \begin{array}{ccc} \oplus^k T\G &\stackrel{\omega}{\to} & \mathbb R \\ \downarrow \downarrow &
 &\downarrow  \\  \oplus^k TM& \stackrel{}{\to} & \{pt\} \\
 \end{array}.$$
 \end{theorem}
 We used the notation $\oplus^k E = E \oplus \ldots \oplus E$, $k$ terms, for $E$ a vector bundle.
The cases $k=1$ and $k=2$ treated in Sections \ref{4.1} and \ref{4.2} are recovered, provided we set $\oplus^0 T\G = \G$ and $\oplus^0 TM = M$. Equations (\ref{1form}) and (\ref{sum2forms}) are particular cases of the following.
\begin{theorem} A $k$-form $\omega$ on $\varGamma$ is {multiplicative} if and only if, for all pairs of composable elements of $T\varGamma$, $(X_1, Y_1)$, $(X_2, Y_2)$, ..., $(X_k, Y_k)$,
\begin{equation}\label{F}
\omega(X_1 \times_{(T\varGamma)} Y_1, X_2 \times_{(T\varGamma)} Y_2,...,X_k \times_{(T\varGamma)} Y_k) = \omega(X_1, X_2 ,...,X_k) + \omega(Y_1, Y_2,...,Y_k).
\end{equation}
\end{theorem}

\subsection{Multiplicative $k$-vectors on Lie groupoids}
Similarly, by definition, a $k$-vector on $\varGamma$ is {\it multiplicative} if it defines
a morphism of Lie groupoids,
$$ \begin{array}{ccc} \oplus^{k-1} T^* \G &\stackrel{}{\to} & T\G \\ \downarrow \downarrow &
 &\downarrow \downarrow  \\  \oplus^{k-1} A^*& \stackrel{}{\to} & TM \\
 \end{array},$$
or, equivalently,
$$ \begin{array}{ccc} \oplus^k T^* \G &\stackrel{}{\to} & \mathbb R \\ \downarrow \downarrow &
 &\downarrow   \\  \oplus^k A^*& \stackrel{}{\to} & \{pt\} \\
 \end{array}.$$
The cases $k=1$ and $k=2$ treated in Sections \ref{4.1} and \ref{3.4} are recovered provided we set $\oplus^0 T^*\G = \G$ and $\oplus^0 A^* = M$. 

Applying a multiplicative $k$-form to a multiplicative $k$-vector yields a multiplicative function.
\section{Linear and morphic multivectors on Lie algebroids}\label{linearmorphic}
In this section we shall study the infinitesimal counterparts of multiplicative multivectors on Lie groupoids, which are multivectors on Lie algebroids that have been called ``morphic''.
To this end, we shall first study the linear vector fields, then the linear multivectors of any degree which can be defined on any vector bundle. Then, we shall assume that the vector bundle under consideration is a Lie algebroid and we shall define and characterize in various ways the morphic 
multivectors. In the next section, we shall study the dual case of the morphic differential forms on Lie algebroids  which are the infinitesimals of multiplicative differential forms on Lie groupoids.
\subsection{Linear vector fields on vector bundles}
A function on a vector bundle is called linear if it is a morphism of vector bundles over $M$, $E \to M \times {\mathbb R}$, i.e., in any local coordinate system, it is linear in the fiber coordinates.

\begin{proposition}Let $X$ be a vector field on a vector bundle, $E \to M$.
The following properties are equivalent:

\noindent(i) $X$ is projectable to a vector field ${\underline X}$ on $M$ and preserves linear functions.

\noindent(ii) $X : E \to TE$ is a vector-bundle morphism over a vector field ${\underline X} : M \to TM$,
$$ 
\begin{array}{ccc}  E &\stackrel{X}{\to} & TE  \\ \downarrow & \downarrow
 &\downarrow  \\ M & \stackrel{{\underline X}}{\to} & TM \\
 \end{array}.
 $$
\vspace{-.3cm}

\noindent(iii) $X$ defines a vector-bundle morphism,
$$ 
\begin{array}{ccc}  T^*E &\stackrel{X}{\to} & {\mathbb R}  \\ \downarrow  &
 &\downarrow  \\ E^* & \stackrel{}{\to} & \{pt\} \\
 \end{array}.
 $$
\vspace{-.3cm}

\noindent(iv) $X$ generates a flow of local vector-bundle automorphisms of $E$ over the 
flow of a vector field ${\underline X}$ on $M$.

\noindent(v) In local coordinates, $u = (x,y) = (x^i, y^\alpha)$ on $E$, $X_u = {\underline X}^i(x)\partial_i + X^\alpha_\beta(x)y^\beta \partial_\alpha$ and ${\underline X}_x = {\underline X}^i(x)\partial_i$.
\end{proposition}

If these properties are satisfied, 
$X$ is called an infinitesimal automorphism of the vector bundle, $E$, or a {\it linear vector field} on $E$ over ${\underline X}$.

\begin{definition} A differential operator, $D$, on the sections of $E$ is called a {\it vector-bundle derivation} if there exists a vector field, ${\underline D}$, on $M$ such that, for all functions $f \in C^\infty (M)$  and all sections $\psi$ of $E$, 
\begin{equation}\label{eqCDO}
D(f \psi) = f D\psi + ({\underline D}f) \psi.
\end{equation} 
\end{definition}
\begin{proposition}
 There is a bijective correspondence between linear vector fields on $E$ and vector-bundle derivations of $E$.
\end{proposition}
\noindent{\it Proof}. A linear vector field, $X$, on $E$ over a vector field, ${\underline X}$, on $M$, defines a vector-bundle derivation, $D_X$, by setting, for $\psi$ a section of $E$ and $x \in M$,  
$$ 
(D_X \psi)_x = (T_x\psi)({\underline X}_x) - X_{\psi(x)}.
$$
The right-hand side is a vector tangent to the fiber of $E$ over $x$, at the point $\psi(x)$, and can be canonically identified with an element of the fiber of $E$ over $x$. Conversely, if $D$ is a vector-bundle derivation of $E$ over a vector field, ${\underline D}$, on $M$, we define a vector field on $E$ as follows. Let $y$ be in the fiber of $E$ at $x$ and let $\psi$ be a section of $E$ such that $\psi(x) = y$. 
Set 
$$
(X_D)_y = (T_x\psi)({\underline D}_x) - (D\psi)_x.
$$
Because $D$ is a vector-bundle derivation, the right-hand side is a vector tangent at $y$ to the fiber of $E$ which depends only on $y$ and not on the choice of $\psi$, and which can be identified with an element of the fiber of $E$ at $x$. The vector field $X_D$ on $E$ defined in this way is a linear vector field  over ${\underline D}$, and the maps $X \mapsto D_X$ and $D \mapsto X_D$ are clearly inverses of one another. \qed

The map $X \mapsto D_X$ is an isomorphism of Lie algebras from the algebra of linear vector fields equipped with the Lie bracket of vector fields to the Lie algebra of vector-bundle derivations equipped with the commutator.

A linear first-order differential operator on  $E$ is the Lie derivation of sections with respect to a linear lifting of the vector field ${\underline X}$ on $M$ if and only if it is a vector-bundle derivation over ${\underline X}$.

\begin{remark} A linear vector field on $E$ also defines a vector-bundle derivation on $E^*$ \cite{MX98}.
\end{remark}

\begin{remark} The symbol of a first-order differential operator, $D$, on $E \to M$ 
is the map, $\sigma(D) : T^*M \to \End(E)$ such that  $\sigma(D)(\d_x f) \psi_x = D(f\psi)(x)$, when $\psi$ is a section of $E$ and $f \in C^\infty(M)$ satisfies $f(x)=0$.
If $X$ is a linear vector field over ${\underline X}$, then the symbol of $D_X$ is the vector-bundle 
morphism defined by $\xi \in T^*M \mapsto \langle {\underline X}, \xi \rangle \, {\mathrm{id}}_E \in \End(E)$. 
Therefore, a linear first-order differential operator on a vector bundle $E$ is a vector-bundle derivation if and only if its symbol takes values in the scalar multiples of the 
identity of $E$  \cite{K76}. 
\end{remark}

In the literature, vector-bundle derivations have been variously called derivative endomorphisms, linear first-order differential operators with scalar symbol (see \cite{K76}, \cite{K80} and references therein), CDO's (acronym for ``covariant differential operators''), derivations (\cite{M87} \cite{MX98}, see \cite{M05}), and more recently ``quasi-derivations'' by Grabowska, Grabowski and Urba\'nski \cite{GGU03}. 
For linear vector fields and linear forms, see also \cite{KU99} and \cite{BC12}, and see \cite{KMS93} citing Kol\'a\v r (1982). The following proposition is from \cite{K76}.

\begin{proposition}\label{degree_1} Vector-bundle derivations on $E$ are in bijective correspondence with derivations of degree $0$ of the algebra of sections of the exterior algebra of $E$, $(\wedge^\bullet E, \wedge)$.
\end{proposition}

\subsection{Linear multivectors on vector bundles}
We shall now consider the more general case of linear multivectors.
By definition, a $k$-vector ($k \geq 1$) on a vector bundle, $E$, is {\it linear} if it defines a morphism of vector bundles,
$$ \begin{array}{ccc} \oplus^{k-1} T^* E &\stackrel{}{\to} & TE \\ \downarrow &
 &\downarrow  \\  \oplus^{k-1} E^*& \stackrel{}{\to} & TM \\
 \end{array},$$
or, equivalently,
$$ \begin{array}{ccc} \oplus^k T^* E &\stackrel{}{\to} & \mathbb R \\ \downarrow &
 &\downarrow  \\  \oplus^k E^*& \stackrel{}{\to} & \{pt\} \\
 \end{array}.$$
In fact, these conditions make sense for $k=1$, provided one sets $\oplus^{0} T^* E = E$ and $\oplus^{0} E^*  = M$, when they coincide with the conditions defining a linear vector field.

\begin{example} A bivector $\pi$ on a vector bundle $E$ (a section of $\wedge^2 TE \to E$)  
is called {\it linear} if it defines a skew-symmetric bracket such that the bracket of linear functions (functions linear on the fibers) is a linear function. This condition is satisfied if and only if $\pi$ defines a morphism of vector bundles from $T^*E$ to $TE$ over a morphism $\underline \pi$ of vector bundles from $E^*$ to $TM$, i.e., if and only if it is a linear $2$-vector.
When $\pi$ is a linear Poisson bivector on $E$, it defines a Lie algebroid structure on $E^*$, and $\underline \pi : E^* \to TM$
is its anchor, $\rho_*$.
\end{example}

We remark that a derivation of degree $k - 1$, $D$, of $\Gamma(\wedge^\bullet E)$ is defined by its restrictions to the spaces of elements of degree 0 and of degree  1, $D_0 : C^\infty(M) \to \Gamma(\wedge^{k- 1} E)$ and 
$D_1 : \Gamma(E) \to \Gamma(\wedge^{k} E)$, which satisfy the relation,
\begin{equation}\label{derivation}
D_1(f \psi) = f D_1(\psi) + (D_0f) \wedge \psi,
\end{equation} 
whence the following theorem due to Iglesias-Ponte, Laurent-Gengoux and Xu \cite{ILX12}.

\begin{theorem}\label{degree_k} Linear $k$-vectors ($k \geq 1$) on $E$ are in bijective correspondence with derivations of degree $k - 1$ of the algebra of sections of the exterior algebra of $E$, $(\wedge^\bullet E, \wedge)$. 
\end{theorem}
 
\begin{remark}
The following terminology was introduced in \cite{ILX12}. An ``{\it almost $k$-differential}'' on $E$ ($k \geq 1$) is a derivation of degree $k - 1$ of the algebra of sections of the exterior algebra of $E$, 
$(\wedge^\bullet E, \wedge)$  so that a $k$-vector gives rise to an ``almost $k$-differential''. (See also Remark \ref{remark} below.) 
\end{remark}

The Schouten--Nijenhuis bracket of linear multivectors is linear.
To the Schouten--Nijenhuis bracket of a linear $k$-vector and a linear $k'$-vector corresponds the commutator of an ``almost $k$-differential'' with an ``almost $k'$-differential'' which is an ``almost $(k+k'-1)$-differential''. 

\begin{example} 
To a linear Poisson bivector on $E$ corresponds an ``almost $2$-differential'' of square 0, i.e., a differential on   
$(\Gamma(\wedge^\bullet E), \wedge)$.
It is the Chevalley-Eilenberg differential, $d_{E^*}$, of the Lie algebroid $E^*$ defined by the linear Poisson bivector on $E$.
\end{example} 

\subsection{Morphic vectors and morphic multivectors on Lie algebroids}
Infinitesimally multiplicative multivectors on Lie algebroids are called ``morphic''.

\begin{definition}
A vector field on a Lie algebroid $A \to M$ is said to be {\it morphic} if it generates a flow of local Lie-algebroid automorphisms of $A$.
\end{definition}

If $A$ is a Lie algebroid with base $M$, then
$TA$ is a Lie algebroid with base $TM$ \cite{C94} \cite{MX94}.
A vector field on $A$ is morphic if and only if it is a Lie-algebroid morphism \cite{HM90}
from $A \to M$ to $TA \to TM$ over a vector field $M \to TM$.

When $A$ is a Lie algebroid, $\Gamma(\wedge^\bullet A)$ equipped with the exterior multiplication and the Schouten--Nijenhuis bracket, $[~,~]$, is a Gerstenhaber algebra.
The vector-bundle derivation defined by a morphic vector field is a derivation of both the exterior multiplication and the Schouten--Nijenhuis bracket of sections of the exterior algebra of $A$, i.e., it is a derivation of the Gerstenhaber algebra of sections of $\wedge^\bullet A$. 
\begin{proposition} Morphic vector fields on $A$ are in bijective correspondence with derivations of degree $0$ of the Gerstenhaber algebra of $A$, 
$(\Gamma(\wedge^\bullet A), \wedge, [~,~])$.
\end{proposition}
The definition of morphic vector fields is given in \cite{MX98} and it was extended to multivectors in \cite{BCO09} \cite{BC12}. 
If $A \to M$ is a Lie algebroid, then $T^*A$ is a Lie algebroid with base $A^*$. 
In fact, since $A$ is a Lie algebroid, $A^*$ is a (linear) Poisson manifold and therefore $T^*A^* \to A^*$
is a Lie algebroid. Applying the canonical map $T^*A^* \to T^*A$ introduced by Tulczyjew (see \cite{MX94}) yields the vector-bundle structure of $T^*A \to  A^*$ as well as its Lie-algebroid structure.

A linear $k$-vector ($k \geq 1$) on $A$ is called {\it morphic} if it defines a Lie-algebroid morphism,
$$ \begin{array}{ccc} \oplus^{k-1} T^* A &\stackrel{}{\to} & TA \\ \downarrow &
 &\downarrow  \\  \oplus^{k-1} A^*& \stackrel{}{\to} & TM \\
 \end{array},$$
or, equivalently,
$$ \begin{array}{ccc} \oplus^k T^*A &\stackrel{}{\to} & \mathbb R \\ \downarrow &
 &\downarrow  \\  \oplus^k A^*& \stackrel{}{\to} & \{pt\} \\
 \end{array}.$$
In fact, these conditions make sense for $k=1$, provided one sets $\oplus^{0} T^* A = A$ and $\oplus^{0} A^*  = M$, when they coincide with the conditions defining a morphic vector field. The following was proved in \cite{ILX12} (see also \cite{LSX11}) and in \cite{BC12}.

\begin{theorem} Morphic $k$-vectors ($k \geq 1$) on $A$ are in bijective correspondence with derivations of degree $k - 1$ of the Gerstenhaber algebra of $A$, $(\Gamma(\wedge^\bullet A), \wedge, [~,~])$.
\end{theorem} 

\begin{remark}\label{remark} A derivation of degree $k - 1$ ($k \geq 1$) of the Gerstenhaber algebra, $(\Gamma(\wedge^\bullet A), \wedge, [~,~])$, is called a ``$k$-{\it differential}'' in \cite{ILX12}. Since it has square 0 only if the associated morphic multivector, $X$, is {\it homological}, i.e., if $[X,X]=0$,
 where $[~,~]$ is the Schouten--Nijenhuis bracket, it would be preferable to speak of ``$k$-derivation'' in general and to reserve the term 
 ``$k$-differential'' for the case of a morphic multivector which is homological. We also remark that, in this terminology, a differential in the usual sense is a ``$2$-differential''.
\end{remark}
 
\begin{example} To a morphic Poisson bivector on $A$ corresponds a derivation of degree $1$, of square~0, i.e., a differential on   
$(\Gamma(\wedge^\bullet A), \wedge, [~,~])$. It is the Chevalley-Eilenberg differential of the Lie algebroid $A^*$ defined by the morphic, hence linear, Poisson bivector on $A$. Since it is a derivation of the Gerstenhaber bracket of $\Gamma(\wedge^\bullet A)$, the pair $(A,A^*)$ is a Lie bialgebroid.  
In fact, $(A,A^*)$ is a Lie bialgebroid if and only if the space of sections of the exterior algebra of $A$ is a differential Gerstenhaber algebra under the Chevalley-Eilenberg differential of $A^*$ \cite{K95} \cite{X99} \cite{ILX12}.
That $(A,A^*)$ is a Lie bialgebroid if and only if the vector-bundle isomorphism $T^*A \to TA$, induced by the Lie-algebroid structure of $A^*$ composed with the canonical map $T^*A^* \to T^*A$, is a morphism of Lie algebroids was proved in \cite{MX94}. 
\end{example}

\smallskip

\subsection{Infinitesimals of multiplicative multivectors on Lie groupoids}
 Applying the Lie functor to a multiplicative multivector, $X$, on a Lie groupoid, $\G$, yields the restriction of the tangent lift of $X$ (see \cite{YI73} \cite{GU95} \cite{GU97}) to the Lie algebroid of $\G$.  To the multiplicativity property of a multivector on a Lie groupoid, there corresponds the fact that its infinitesimal counterpart is a derivation of the bracket of its Lie algebroid. 
\begin{theorem}\label{infinitesimal} The infinitesimal of a multiplicative $k$-vector ($k \geq 1$) on a Lie groupoid, $\G$, is a morphic $k$-vector on the Lie algebroid, $A$, of $\G$.
If the source map of $\G$ has simply-connected fibers, each morphic $k$-vector on $A$ is the infinitesimal of a multiplicative $k$-vector on $\G$.
\end{theorem}
\noindent{\it Idea of proof}. To determine the morphic $k$-vector on $A$ associated to a given multiplicative $k$-vector on $\G$, one determines an associated $k$-differential on $A$. Recall that the Gerstenhaber bracket of a multiplicative multivector and a right-invariant multivector on a Lie groupoid is right-invariant \cite{MX94}. Therefore, given a multiplicative $k$-vector, $X$, on $\G$, to a section $u$ of $\wedge^\ell A$, prolonged into a right-invariant section $u^\rho$ of $\wedge^\ell T\G$, we associate the section $\d_X u$ of $\wedge^{\ell + k - 1} A$ whose right-invariant prolongation is the Schouten--Nijenhuis bracket $[X, u^\rho]$. Thus, to any multiplicative $k$-vector, $X$, on $\G$, there corresponds an endomorphism $\d_X$ of degree $k-1$ of $\Gamma(\wedge^\bullet A)$ which is a derivation of the exterior multiplication, and which in turn defines a linear
 $k$-vector, Lie$(X)$, on $A$. The definitions of ``multiplicative'' and ``morphic'' in terms of morphisms show that if $X$ is multiplicative, then 
 Lie$(X)$ is morphic.
The converse property, integrating a given morphic vector field on $A$ into a multiplicative multivector on $\G$, when $\G$ is source-simply-connected, was more involved. It was achieved by Iglesias-Ponte, Laurent-Gengoux and Xu \cite{ILX12}, using $A$-paths. However, interpreting multiplicative multivectors on Lie groupoids as morphisms of Lie groupoids and morphic multivectors as morphisms of Lie algebroids shows that this theorem is a corollary of the theorem on the integration of a Lie-algebroid morphism to a Lie-groupoid morphism in the case of a source-simply-connected Lie groupoid \cite{MX00}.

\begin{example} In the case of a Lie group $G$, to a multiplicative $k$-vector, $X$, on $G$ is associated its intrinsic derivative, $\d_eX : \g \to \wedge^k \g$, which can be extended uniquely to a derivation of degree $k-1$ of $\wedge^\bullet \g$. If $k = 2$, denoting $X$ by $\pi$, $\d_e\pi$ is a cobracket on $\g$ which we denote by $\gamma$ and which can be viewed as a linear bivector field on $\g$. The associated derivation of degree 1 of $\wedge^\bullet \g$ is $\{\gamma, \cdot \}$, where $\{~,~\}$ is the big bracket. If $[\pi,\pi]=0$, the cobracket $\gamma$ defines a Lie bracket on $\g^*$ and the derivation $\{\gamma, \cdot \}$ is the Chevalley-Eilenberg differential of the Lie algebra $\g^*$ (see, e.g., \cite{K92}).
\end{example}

\begin{example} If $\G$ is a Poisson groupoid with Lie algebroid $A$, the infinitesimal of the Poisson bivector of $\G$ is a morphic Poisson bivector on $A$, which defines a Lie algebroid structure on $A^*$ (since it is linear), and the pair $(A,A^*)$ is a Lie bialgebroid (since it is morphic).   
Conversely, under the assumptions of Theorem \ref{infinitesimal}, a Lie bialgebroid can be integrated to a Poisson groupoid, results first proved in \cite{MX94} and \cite{MX98}.
\end{example}

\section{Infinitesimal multiplicativity of forms on Lie algebroids}\label{IM}
We shall now define the linear differential forms on vector bundles -- a particular case of the polynomial differential forms defined in \cite{TU00} -- and the morphic differential forms on Lie algebroids, and we shall then describe the correspondence between multiplicative forms on a Lie groupoid and morphic forms on its Lie algebroid.

\subsection{Linear 1-forms on vector bundles}

\begin{proposition}Let $\xi$ be a 1-form on a vector bundle $E \to M$.
The following properties are equivalent.

\noindent (i) For all linear vector fields, $X$, on $E$, $\langle  \xi, X \rangle $ is a linear function on $E$, i.e. a vector-bundle morphism from $E \to M$ to $M \times {\mathbb R} \to M$.

\noindent (ii) $\xi$ is a vector-bundle morphism $E \to T^*E$ over a section of $E^*$, ${\underline \xi} : M \to E^*$,
$$ \begin{array}{ccc}  E &\stackrel{\xi}{\to} & {T^*E}  \\ \downarrow &
 &\downarrow  \\ M & \stackrel{{\underline \xi}}{\to} & E^* \\
 \end{array}.$$

\noindent (iii)  $\xi : E \to T^*E$ defines a vector-bundle morphism,
$$
 \begin{array}{ccc}  TE &\stackrel{\xi}{\to} & {\mathbb R}  \\ \downarrow &
 &\downarrow  \\ TM & \stackrel{}{\to} & \{pt\} \\
 \end{array}.
 $$

\noindent (iv) In local coordinates $ u = (x,y) =(x^i, y^\alpha)$ on $E$, $\xi_u = \xi_{i\alpha}(x)y^\alpha \d x^i + {\underline \xi}_\alpha(x) \d y^\alpha$, 
and ${\underline \xi}_x = {\underline \xi}_\alpha(x) \d y^\alpha$.
\end{proposition}

If these properties are satisfied, 
$\xi$ is said to be a {\it linear 1-form} over ${\underline \xi}$  \cite{MX98} \cite{KU99}.

\medskip

In local coordinates, the evaluation of a linear $1$-form, $\xi_u= \xi_{i\alpha}(x)y^\alpha dx^i + {\underline \xi}_\alpha(x) dy^\alpha$, on a linear vector field,  $X_u = {\underline X}^i(x)\partial_i + X^\alpha_\beta(x)y^\beta \partial_\alpha$, is the linear function whose value at $u\in E$ is  
$\langle \xi_u, X_u \rangle =  (\xi_{i\alpha}  {\underline X}^i + {\underline \xi}_\beta X^\beta_\alpha) y^\alpha$.

The differential, $\d f$, of a linear function $f : E \to {\mathbb R}$ is a linear 1-form since, for each linear vector field, $X$, on $E$, $\langle  \d f, X \rangle = X \cdot f$ is a linear function.

\medskip
\subsection{Linear $k$-forms on vector bundles}\label{linear}
We shall now consider the more general case of linear forms of higher degree, $k \geq1$.
Here and in Sections \ref{7.3} and \ref{IMform} we follow the articles by Bursztyn, Cabrera and Ortiz, \cite{BCO09} and \cite{BC12}.
The case of 2-forms was investigated before the general case because of its importance in symplectic and presymplecic geometry.
\begin{theorem}Let $\xi$ be a $k$-form ($k \geq 1$) on a vector bundle $E \to M$.
The following properties are equivalent.

\noindent (i) For all linear vector fields, $X_1$, \ldots, $X_k$ on $E$, $\xi (X_1, \ldots, X_k)$ is a linear function.

\noindent (ii) $\xi : E \to \wedge^k T^*E$ defines a vector-bundle morphism, $\oplus^{k-1}TE \to T^*E$, over a vector-bundle morphism, ${\underline \xi} : \oplus^{k-1}TM \to E^*$,
$$ \begin{array}{ccc}  \oplus^{k-1}TE &\stackrel{\xi}{\to} & {T^*E}  \\ \downarrow &
 &\downarrow  \\  \oplus^{k-1}TM & \stackrel{{\underline \xi}}{\to} & E^* \\
 \end{array}.$$

\noindent (iii)  $\xi : E \to \wedge^k T^*E$ defines a vector-bundle morphism,
$$
 \begin{array}{ccc}  \oplus^{k}TE &\stackrel{\xi}{\to} & {\mathbb R}  \\ \downarrow &
 &\downarrow  \\ \oplus^{k}TM & \stackrel{}{\to} & \{pt\} \\
 \end{array}.
$$

\noindent (iv) In local coordinates $u = (x,y) = (x^i, y^\alpha)$ on $E$, 
$$
\xi(x, y) = \frac{1}{k!} \xi_{i_1 \ldots i_k \alpha}(x)y^\alpha \d x^{i_1} \wedge \ldots \wedge \d x^{i_k} + \frac{1}{(k-1)!}{\underline \xi}_{i_1 \ldots i_{k-1} \alpha}(x) \d x^{i_1} \wedge \ldots \wedge \d x^{i_{k-1}} \wedge \d y^\alpha
$$
and ${\underline \xi}_x(\partial_{i_1}, \ldots, \partial_{i_{k - 1}}) = {\underline \xi}_{i_1 \ldots i_{k-1} \alpha}(x) \d y^\alpha$.
\end{theorem} 
If these properties are satisfied, 
$\xi$ is called a {\it linear $k$-form} over ${\underline \xi} : \oplus^{k-1}TM \to E^*$.
In fact, these conditions make sense for $\, k=1$, provided one sets $\oplus^{0}TE = E$ and $ \oplus^{0}TM = M$, when they coincide with the conditions defining a linear $1$-form.

The differential of a linear $k$-form is a linear $(k+1)$-form. 

\begin{example} A linear $2$-form on a vector bundle $E \to M$ is a pair of vector-bundle morphisms, $(\xi,{\underline \xi})$, in the commutative diagram, 
$\begin{array}{ccc} TE  &\stackrel{\xi}{\to} & T^*E  \\ \downarrow &
 &\downarrow  \\ TM & \stackrel{{\underline \xi}}{\to} & E^* \\
 \end{array}$. Equivalently, the graph of $\xi$ is a vector sub-bundle of the vector bundle $TE \oplus T^*E \to TM \oplus E^*$.
  
A linear $2$-form may be considered as a morphism of double vector bundles \cite{M05} from 
$\begin{array}{ccc} TE  &\stackrel{}{\to} & TM  \\ \downarrow &
 &\downarrow  \\ E & \stackrel{}{\to} & M \\
 \end{array}$ to $\begin{array}{ccc} T^*E  &\stackrel{}{\to} & E^*  \\ \downarrow &
 &\downarrow  \\ E & \stackrel{}{\to} & M \\  \end{array}$.
\end{example}
\subsection{The structure of linear $k$-forms}\label{7.3}
We shall determine the structure of linear $k$-forms, $k \geq 1$, using local coordinates. 
For $h \geq 1$, let $\zeta : E \to \wedge^h T^*M$ be a morphism of  vector bundles over $M$. Then, for $u = (x,y) = (x^i,y^\alpha) \in E$, $\zeta(u)$ is an $h$-form on $M$ at $x$, with components $\frac{1}{h!}\zeta_{i_1 \ldots i_{h} \alpha}(x) y^\alpha$. The pull-back of $\zeta(u)$ by the projection $q : E \to M$ is a linear $h$-form on $E$ at $u$, 
$$
\hat{\zeta}(u) = \frac{1}{h!}
\zeta_{i_1 \ldots i_{h} \alpha}(x) y^\alpha \d x^{i_1} \wedge \ldots \wedge \d x^{i_h}.
$$
The differential of $\hat{\zeta}$ at $u$ satisfies 

\noindent{\small{
$(-1)^{h} \d\hat{\zeta}(u)=\frac{1}{(h + 1)!}\partial_{i_{h+1}}\zeta_{i_1 \ldots i_{h} \alpha}(x) y^\alpha \d x^{i_1} \wedge \ldots \wedge \d x^{i_{h+1}}
+ \frac{1}{h!}\zeta_{i_1 \ldots i_{h} \alpha}(x) \d x^{i_1} \wedge \ldots \wedge \d x^{i_{h}} \wedge \d y^\alpha$.}}

\noindent Let $\xi$ be a linear $k$-form on $E$ over ${\underline \xi}$. If $k > 1$, ${\underline \xi}$ is a vector-bundle morphism, from $\wedge 
^{k-1}TM$ to $E^*$, and we denote its transpose by the same symbol. If $k =1$, ${\underline \xi}$ is a section of $E^*$ and we denote the corresponding linear $0$-form on $E$ by the same symbol. We see that 
$\xi - (-1)^{k-1} \d\hat{\underline \xi}$ is a linear $k$-form on $E$ over the zero map, which is itself of the 
form $\hat\nu$ for some (unique) map $\nu : E \to \wedge^{k}T^*M$. Thus 
$$
\xi = \d\hat \mu + \hat\nu,
$$
where 
$$
\mu = (-1)^{k-1}{\underline \xi}.
$$
Therefore, to $(\xi,{\underline \xi})$ is associated a pair of maps $\mu = (-1)^{k-1}{\underline \xi} : E \to \wedge^{k-1}T^*M$ and $\nu : E \to\wedge^{k}T^*M$ such that $\hat \nu = \xi - (-1)^{k-1} \d \hat {\underline \xi}$. Conversely, the pair $(\mu,\nu)$ determines $(\xi,{\underline \xi})$ by $\xi = \d\hat \mu + \hat \nu$.

If, in particular, $\xi$ is closed, $\d\hat\nu = 0$, whence $\nu = 0$ and
$$
\xi =  \d\hat{\mu}.
$$ 
Therefore, in the case of closed linear forms, the map $\underline \xi$ determines the form $\xi$.

\begin{remark} For any vector-bundle morphism, $\zeta : E \to T^*M$, the $1$-form $\hat \zeta$ on $E$ is the pull-back of the Liouville $1$-form $\theta = p_i \d x^i$ of $T^*M$ by the map $\zeta$, 
$$
\hat\zeta = \zeta^*\theta,
$$
as can be seen in local coordinates. If $(\xi,{\underline \xi})$ is a linear $2$-form and if $\xi$ is closed, $\xi = - \d\hat{\underline \xi}$. In that case, since  
$\hat{\underline \xi} = {\underline \xi}^*\theta$, we obtain  $\xi = - \d\hat{\underline \xi} = - \d({ \underline \xi}^*\theta)= - {\underline \xi}^*(\d\theta)$.  Therefore,
$$
\xi = {{\underline \xi}}^* \omega,
$$ 
where $\omega = - \d\theta = \d x^i \wedge \d p_i$ is the canonical symplectic $2$-form on $T^*M$ (see \cite{KU99}).
\end{remark}
\subsection{Morphic forms on Lie algebroids}\label{IMform}
When the vector bundle under consideration is a Lie algebroid, $A$, one can define the {\it morphic} $1$-forms as those linear 1-forms, $(\xi, {\underline \xi})$, which are such that the following diagram is a morphism of Lie algebroids, 
$$ \begin{array}{ccc}  A &\stackrel{\xi}{\to} & {T^*A}  \\ \downarrow &
&\downarrow  \\ M & \stackrel{{\underline \xi}}{\to} & A^* \\
\end{array},$$
or, equivalently,
$$
 \begin{array}{ccc}  TA &\stackrel{\xi}{\to} & {\mathbb R}  \\ \downarrow &
&\downarrow  \\ TM & \stackrel{}{\to} & \{pt\} \\
\end{array}
$$
is a morphism of Lie algebroids.

More generally, the {\it morphic} $k$-forms are those linear $k$-forms, $(\xi, {\underline \xi})$, which are such that the following diagram is a morphism of Lie algebroids, 
$$ \begin{array}{ccc}  \oplus^{k-1}TA &\stackrel{\xi}{\to} & {T^*A}  \\ \downarrow &
&\downarrow  \\  \oplus^{k-1}TM & \stackrel{{\underline \xi}}{\to} & A^* \\
\end{array},$$
or, equivalently,
$$
 \begin{array}{ccc}  \oplus^{k}TA &\stackrel{\xi}{\to} & {\mathbb R}  \\ \downarrow &
&\downarrow  \\ \oplus^{k}TM & \stackrel{}{\to} & \{pt\} \\
\end{array}
$$
is a morphism of Lie algebroids.

\begin{example} A morphic 2-form on $A$ is a pair, $(\xi,\underline \xi)$, of vector-bundle morphisms, $\xi : TA \to T^*A$ and $\underline \xi : TM \to A^*$, defining a Lie-algebroid morphism, or a pair, $(\mu,\nu)$, $\mu : A \to T^*M$ and $\nu : A \to \wedge^2 T^*M $. The map $\mu$ is identified with the 
transpose of  $- \underline \xi$. 
\end{example}
 
We have seen that a morphic $k$-form,  $(\xi,{\underline \xi})$, is defined by a pair $(\mu, \nu)$ of maps from $A$ to $\wedge^{k-1} T^*M$ and from $A$ to 
$\wedge^{k} T^*M$, respectively. The question is to translate the Lie-algebroid morphism property of the pair of maps, $(\xi, {\underline \xi})$, $\xi : 
\wedge^{k-1}TA \to T^*A$ and ${\underline \xi} : \wedge^{k-1}TM \to A^*$, into equivalent properties for the pair of maps $(\mu,\nu)$.
Recall that $\mu$ is equal, up to sign, to ${\underline \xi}$, when ${\underline \xi}$ is identified with its transpose. 
\begin{theorem}\label{BC2012} Let $A$ be a Lie algebroid over $M$, with bracket $[~,~]$ and anchor $\rho$, and let $(\xi, {\underline \xi})$ be a linear $k$-form ($k \geq 1$) on $A$. Let $\mu = (-1)^{k-1} {\underline \xi}$ and let $\xi  = \d\hat \mu + \hat \nu$. The pair $(\xi, {\underline \xi})$ is a morphic $k$-form  on  $A$ if and only if $\mu : A \to \wedge^{k-1} T^*M$ and $\nu : A \to \wedge^{k} T^*M$
satisfy the relations, for all sections $u$ and $v$ of $A$,

\noindent(i) $i_{\rho (u)}(\mu(v)) = - i_{\rho (v)}(\mu(u))$,

\noindent(ii) $\mu[u,v] = L_{\rho (u)}(\mu v) - i_{\rho (v)} \d(\mu u) - i_{\rho (v)}(\nu u)$,

\noindent(iii) $\nu[u,v]= L_{\rho (u)}(\nu v) - i_{\rho (v)} \d(\nu u)$.
\end{theorem}

The pair $(\mu, \nu)$ is what is called in the literature an ``IM form'' (short for infinitesimally multiplicative). 

\begin{remark} In view of Theorem \ref{morphic} below, it would be more logical to call the pair $(\xi, {\underline \xi})$, instead of 
the pair $(\mu, \nu)$, an infinitesimally multiplicative form.
\end{remark}

Recall that if the morphic $k$-form, $\xi$, is closed, then $\nu = 0$ and therefore $\xi$ is entirely defined by the map $\mu : A \to \wedge^{k-1} T^*M$.  In that case, it is the map $\mu$ itself which is called an IM form. The above equations simplify and are reduced to two.
\begin{corollary} 
Closed morphic $k$-forms on a Lie algebroid $A \to M$ are in bijective correspondence with maps $\mu : A \to \wedge^{k-1} T^*M$ such that,  for all sections $u$ and $v$ of~$A$,

\noindent(i) $i_{\rho (u)}(\mu v) = - i_{\rho (v)}(\mu u)$,

\noindent(ii) $\mu[u,v] = L_{\rho (u)}(\mu v) - i_{\rho (v)} \d(\mu u)$.
\end{corollary}

A linear $k$-form, $\xi$, on a Lie algebroid, $A \to M$, with anchor $\rho$, is called relatively $\phi$-closed, where $\phi$ is a closed $(k+1)$-form on $M$, if the $(k+1)$-form $\d \xi$ evaluated at $u \in A$ is the differential of the pull-back to $A$ of the $k$-form on $M$, $i_{\rho (u)}\phi$.
For the case of closed and relatively closed $2$-forms, see~\cite{BCWZ04}. 
The proof of Theorem \ref{BC2012} for arbitrary $k \geq 1$ was given by Bursztyn and Cabrera in~\cite{BC12}.

\subsection{Infinitesimals of multiplicative forms on Lie groupoids}
What is the infinitesimal of a multiplicative $k$-form on a Lie groupoid $\varGamma$?
Applying the Lie functor to a $k$-form $\omega$, one obtains a $k$-form Lie$(\omega)$ on the Lie algebroid, $A$, of $\G$.  
This $k$-form on $A$ is linear in the sense made precise in Section \ref{linear}.
What are its properties when $\omega$ is multiplicative?
To the multiplicativity property of a form on a Lie groupoid, there corresponds the fact that its infinitesimal counterpart is an infinitesimally multiplicative form on its Lie algebroid.

\begin{theorem}\label{morphic} The infinitesimal of a multiplicative $k$-form on a Lie groupoid $\G$ is a morphic $k$-form on the Lie algebroid, $A$, of $\G$.
If the source map of $\G$ has simply-connected fibers, each morphic $k$-form on $A$ is the infinitesimal of a multiplicative $k$-form on $\G$.
\end{theorem}
\noindent{\it Idea of proof}. 
A multiplicative $k$-form, $\omega$, on $\G$ induces a linear $k$-form, ${\mathrm{Lie}}\,( \omega)$, on $A$ which is the restriction to $A$ of the tangent lift of $\omega$ (see \cite{GU95} \cite{GU97}). The definitions of ``multiplicative'' and ``morphic'' in terms of morphisms show that if $\omega$ is multiplicative, then ${\mathrm{Lie}}\,( \omega)$ is morphic.
The converse property, integrating a given morphic form on $A$ into a multiplicative form on $\G$, when $\G$ is source-simply-connected, is more involved. However, interpreting multiplicative forms on Lie groupoids as morphisms of groupoids, and morphic forms as morphisms of Lie algebroids yields a proof of the converse as a corollary of the theorem on the integration of a Lie-algebroid morphism to a Lie-groupoid morphism in the case of a source-simply-connected Lie groupoid \cite{MX00}.

Let $\omega$ be a multiplicative $2$-form on $\varGamma$ and let $(\mu,\nu)$ be the IM form associated to $(\xi, \underline \xi) = {\mathrm{Lie}}\, (\omega)$. Various particular cases arise.
If $\omega$ is closed, $\nu = 0$ (the presymplectic case).
If $\omega$ is relatively $\phi$-closed, $\nu$ is the map $u \in A \mapsto i_{\rho (u)} \phi \in  \wedge^2 T^*M$.
If $\omega$ is non-degenerate, $\underline \xi : A \to T^*M$ is an isomorphism.  
If $\omega$ is non-degenerate and closed (the symplectic case), the isomorphism $\underline \xi : A \to T^*M$ determines a Poisson structure on $M$ such that the source map is a Poisson map and the target map is anti-Poisson.
If $\omega$ is non-degenerate and relatively $\phi$-closed, the isomorphism $\underline \xi : A \to T^*M$ determines a ``$\phi$-twisted Poisson structure'' on $M$.

\subsection{Lie groupoid and Lie algebroid cocycles}
Once the double complex of a Lie groupoid and the cohomology of a Lie algebroid are defined, multiplicative forms and multiplicative multivectors on Lie groupoids can be viewed as Lie groupoid cocycles, while morphic multivectors and morphic forms on Lie algebroids are Lie algebroid cocycles. 
Crainic \cite{C04} observed that if a multiplicative $2$-form, $\omega$, on $\G$ is closed, then Lie$(\omega)$ is a Lie algebroid cocycle and Lie$(\omega)(u,v) = \langle \mu (u), \rho (v) \rangle$, for all sections $u$ and $v$ of $A$.
Under assumptions of connectedness, there is a general ``Van Est theorem'' relating the cohomology of a line of the double complex of a Lie groupoid to that of a line in the bigraded Weil algebra of its Lie algebroid. The derivation and integration theorems summarized above then appear as particular cases of this general theorem. See Camillo Arias Abad and Crainic 
\cite{AC11} and Jotz \cite{J14}. These results are efficiently formulated in terms of supermanifolds \cite{R15}.

\subsection{Bisections}
A {\it bisection} of a Lie groupoid $\varGamma \overrightarrow{\rightarrow} M$ is a section, $M \to     \varGamma$, of the source map $\alpha$ such that its composition with the target map $\beta$ is a diffeomorphism of the base. Bisections act on tensor fields and differential forms on the groupoid by the associated left and right translations. 
The translations defined by the bisections in a Lie groupoid play the role of the translations on a Lie group defined by an element in the group.
By definition, each bisection, $\Sigma$, of $\varGamma$ induces a diffeomorphism $\phi_\Sigma = \beta \circ \Sigma$ of the base, $M$, of $\G$. 
Therefore a bisection, $\Sigma$, acts on $\Omega^k (M)$, the $k$-forms on $M$, by 
$(\Sigma, \lambda) \mapsto \Sigma \cdot \lambda = \phi_\Sigma^* (\lambda)$, for $\lambda \in \Omega^k (M)$.
In other words, for vector fields, $\varepsilon_1 , \varepsilon_2,...  , \varepsilon_k $, tangent to $M$,
$$(\Sigma \cdot \lambda)(\varepsilon_1 , \varepsilon_2  ,...,  \varepsilon_k) = \lambda(T\phi_\Sigma (\varepsilon_1) , T\phi_\Sigma (\varepsilon_2),... ,T\phi_\Sigma (\varepsilon_k)).
$$

The {bisections} of $\varGamma$ form a group, denoted by ${\mathcal G}(\varGamma)$, with multiplication defined by:
$$
(\Sigma_1 * \Sigma_2) (x) = \Sigma_1(\beta \Sigma_2 (x)) \times_{(\varGamma)} \Sigma_2 (x),
$$ 
for $\Sigma_1, \Sigma_2 \in {\mathcal G}(\varGamma)$, and $x \in M$. Here $\times_{(\varGamma)}$ denotes the groupoid multiplication.
Since $\phi_{\Sigma_1 * \Sigma_2} = \phi_{\Sigma_1} \circ \phi_{\Sigma_2}$, the action thus defined is an action on $\Omega^k (M)$ of the opposite group ${\mathcal G}(\varGamma)^{opp}$ of ${\mathcal G}(\varGamma)$.

Using the concept of bisection of a Lie groupoid, Xu \cite{X95} proved the following theorem showing that a bivector satisfying the morphism property satisfies a relation that generalizes a property of affine bivectors on Lie groups.
\begin{theorem} 
Let $\pi$ be a bivector on a Lie groupoid. If $\pi$ satisfies the morphism property, then 
$$\pi_{gh} = \tilde g \cdot \pi_h + \pi_g \cdot \tilde h -  \tilde g \cdot \pi_{\alpha g}  \cdot \tilde h,$$
for all $g, h \in     \varGamma$ such that $\alpha g= \beta h$, where $\tilde g$ is a bisection that takes the value $g$ at $\alpha g$, and $\tilde h$ is a bisection that takes the value 
$h$ at $\alpha h$. 
\end{theorem}

\subsection{Multiplicative forms induce cocycles on the group of bisections}
We again denote the source map of $\varGamma$ by $\alpha$. For $\Sigma \in {\mathcal G}(\varGamma)$, $\alpha_{|\Sigma}$ is a diffeomorphism from $\Sigma$ to $M$.
Let $\omega \in \Omega^k (\varGamma)$ be a $k$-form on $\varGamma$.
We define a map, 
$$
c_\omega :  \Sigma \in {\mathcal G}(\varGamma)^{opp} \mapsto (\alpha_{|\Sigma}^{-1})^* \omega \in \Omega^k (M).
$$
The following unpublished result was communicated to me by Laurent-Gengoux \cite{L15}.
\begin{theorem} 
If $\omega$ is a multiplicative $k$-form on $\varGamma$, then $c_\omega$ is a {1-cocycle}
on the group $ {\mathcal G}(\varGamma)^{opp}$ with values in the ${\mathcal G}(\varGamma)^{opp}$-module $\Omega^k (M)$,
$$
c_\omega(\Sigma_1 * \Sigma_2) = c_\omega(\Sigma_2) + \Sigma_2 . c_\omega(\Sigma_1).
$$
\end{theorem}

\noindent{\it Proof.}
Let $ \epsilon$ be a vector field on $M$. 
To evaluate $T(\alpha_{|\Sigma_1 *  \Sigma_2}^{-1})\epsilon$, for $\epsilon$ tangent to $M$, we write $\epsilon$ as the tangent at $t=0$ of a curve $x_t$ on $M$. 
Then the image of $x_t$ under $\alpha_{|\Sigma_1 *  \Sigma_2}^{-1}$ is $(\alpha_{|\Sigma_1}^{-1}(\beta \Sigma_2 x_t) \times_{(\varGamma)} (\alpha_{|\Sigma_2}^{-1} x_t)$.
The image of $\epsilon$ under $T(\alpha_{|\Sigma_1 *  \Sigma_2}^{-1})$ is the tangent at $t=0$ of this curve on $\varGamma$, 
 defined by the product in $T\varGamma$,
$$
T(\alpha_{|{\Sigma_1 *  \Sigma_2}}^{-1})\epsilon
 =  \frac{\d}{\d t}_{|t=0}\left((\alpha_{|\Sigma_1}^{-1}(\beta \Sigma_2 x_t) 
 \times_{(\varGamma)} (\alpha_{|\Sigma_2}^{-1} x_t\right) 
$$
$$
 = T(\alpha_{|\Sigma_1}^{-1}) T\varphi_{\Sigma_2} \epsilon \,  \times_{(T\varGamma)} \, T(\alpha_{|\Sigma_2}^{-1}) \epsilon, 
$$  
where $\varphi_{\Sigma_2} = \beta \circ \Sigma_2$ is the diffeomorphism of $M$ defined by $\Sigma_2$.

Now, let us assume that $\omega$ is multiplicative.
Let $ \epsilon_1, \epsilon_2,..., \epsilon_k$ be vector fields on $M$. Apply formula (\ref{F}), with $X_i = T(\alpha_{|\Sigma_1}^{-1})T\varphi_{\Sigma_2} \epsilon_i$ and $Y_i = 
T(\alpha_{|\Sigma_2}^{-1}) \epsilon_i$, $i = 1, \ldots, k$, to obtain, 
$$c_\omega(\Sigma_1 * \Sigma_2) (\epsilon_1, \epsilon_2,..., \epsilon_k) 
$$
$$
= \omega( T(\alpha_{|\Sigma_1}^{-1})T\varphi_{\Sigma_2}\epsilon_1, T(\alpha_{|\Sigma_1}^{-1})T\varphi_{\Sigma_2}\epsilon_2,..., T(\alpha_{|\Sigma_1}^{-1})T\varphi_{\Sigma_2}\epsilon_k)
$$
$$
+ \, \omega (T(\alpha_{|\Sigma_2}^{-1})\epsilon_1, T(\alpha_{|\Sigma_2}^{-1})\epsilon_2,...,T(\alpha_{|\Sigma_2}^{-1})\epsilon_k). 
$$
The first term of the right-hand side is $\Sigma_2 . c_\omega({\Sigma_1}) (\epsilon_1, \epsilon_2, ..., \epsilon_k)$, and the second term is
$c_\omega({\Sigma_2}) (\epsilon_1, \epsilon_2, ..., \epsilon_k)$. \qed


\section{Generalized geometry}\label{GG}
As Nigel Hitchin has shown~\cite{H03}, 
the theory of {generalized complex structures} on the {generalized tangent bundle} of a manifold unifies the symplectic and the complex geometry of the manifold. In this section, we shall examine the case where the base manifold is a Lie groupoid.

\subsection{Generalized tangent bundle} 
The  {\it generalized tangent bundle} of a manifold was introduced by Theodore Courant and Weinstein \cite{CW88} and it was the main subject of Courant's thesis~\cite{C90}. For an arbitrary manifold, $M$, there is a bracket on the sections of $TM \oplus T^*M$, the {\it Courant bracket}, which is skew-symmetric but does {not} satisfy the Jacobi identity.  
The vector bundle $TM \oplus T^*M$ is called the {\it generalized tangent bundle} of $M$, and is also called the {\it standard Courant algebroid}, or the {\it Pontryagin bundle} of $M$,
because of its r\^ole in control theory (see \cite{YM06}).

Consider the following bracket on the sections of the vector bundle $ {\mathbb T} M = TM \oplus T^*M$ over $M$, defined by
$$
{[X+\xi, Y+\eta]= [X,Y] + {\mathcal L}_X \eta - i_Y d\xi},
$$
for $X,Y\in C^\infty(TM)$, $\xi,\eta\in C^\infty(T^*M)$.
In particular, $[X,\xi] =  {\mathcal L}_X \xi$ and $[\xi,X] =  - i_X  d\xi$.
It satisfies the Jacobi identity in its Leibniz form, 
\begin{equation}\label{leibniz}
[u,[v,w]] = [[u,v],w]+[v,[u,w]],
\end{equation}
for all sections $u,v,w$ of the vector bundle ${\mathbb T} M \to M$, but it is not skew-symmetric unless $M$ is a point.
It is called the \textit{Dorfman bracket}\footnote{This bracket is named after Irene Dorfman (1948-1994). It
had appeared in her articles on integrable systems \cite{D87} and in her book on Dirac structures \cite{D93}.
It was introduced in the theory of Courant algebroids circa 1998, by Pavol \v Severa, Xu, and myself, independently, and published only much later. On the development of the theory of Courant algebroids, see \cite{K13}.}. A posteriori, the Courant bracket has been recognized as the skew-symmetrization of the Dorfman bracket.

\subsection{Courant algebroids}
In 1997, Zhang-Ju Liu, Weinstein and Xu defined a Courant bracket on the double of a Lie bialgebroid, thus generalizing the concept of generalized tangent bundle, and they further defined a general concept of {Courant algebroid} \cite{LWX97}. 
The definition of Courant algebroids can be conveniently expressed in terms of {\it Leibniz algebras},
a concept due to Jean-Louis Loday (1946-2012). (Leibniz algebras are also called {\it Loday algebras} or {\it Loday-Leibniz algebras}.)

A {\it  Leibniz algebra} \cite{L93} is a vector space, $L$, equipped with a bilinear bracket, $[~,~] : L \times L \to L$,  such that, for all $u \in L$, $[u,\cdot \,]$ is a derivation of $(L, [~,~])$,
i.e., equation (\ref{leibniz}) is satisfied for all $u, v, w \in L$. This bracket is not, in general, skew-symmetric.

A {\it Courant algebroid} is a vector bundle, $E \to M$, with a fiberwise metric (a non-degenerate symmetric, fiberwise bilinear form), denoted by $\langle ~,~ \rangle$, whose vector space of sections is a Leibniz algebra, with a bracket denoted by $[~,~]$, together with a vector-bundle morphism called the {\it anchor}, $\rho : E \to TM$, such that, for all sections $u$, $v$, $w$ of $E$,

(i) $\rho (u) \langle v,w \rangle = \langle u, [v,w] + [w,v] \rangle$,

(ii) $\rho (u)  \langle v,w \rangle =   \langle [u, v],w \rangle +  \langle v, [u,w] \rangle .$

\noindent The bracket on the space of sections is called the {\it Courant-Dorfman bracket} or, simply, the {\it Dorfman bracket}, and the skew-symmetrization of the 
Dorfman bracket is called the {\it Courant bracket}.

\subsection{Dirac structures}
Let $E \to M$ be a Courant algebroid. A vector sub-bundle of $E$ on which the metric vanishes is called 
{\it isotropic}.\footnote{The term ``isotropic'' is used when dealing with either symmetric or skew-symmetric bilinear forms.} A {\it Dirac bundle} in $E$ is a maximally isotropic sub-bundle of $E$ whose space of sections is closed under the Dorfman bracket (or, equivalently, the Courant bracket). The anchor and the Dorfman bracket of $E$ induce a Lie-algebroid structure on any Dirac bundle in $E$. (Therefore a Dirac bundle is also sometimes called a ``Dirac algebroid''.)

In the case of the standard Courant algebroid, a Dirac bundle in ${\mathbb T}M = TM \oplus T^*M$ is said to define a {\it Dirac structure} on $M$, and $M$ itself is then called a {\it Dirac manifold}.

\begin{example} 
Graphs of bivectors on $M$ are the maximally isotropic sub-bundles, $L$, of $TM \oplus T^*M$ such that $L \cap TM = \{ 0 \}$.
A bivector, $\pi$, on $M$ is Poisson if and only if the graph of $\pi^\sharp$ defines a Dirac structure on $M$. 
\end{example}
\begin{example} 
Graphs of 2-forms on $M$ are the maximally isotropic sub-bundles, $L$, of $TM \oplus T^*M$ such that $L \cap T^*M =  \{ 0 \}$.
A $2$-form, $\omega$, on $M$ is closed if and only if the graph of $\omega^\flat$ defines a Dirac structure on $M$. 
\end{example} 
\subsection{Dirac groupoids}\label{ortiz}
Dirac groupoids, which generalize both Poisson and presymplectic groupoids, were defined by Ortiz in his IMPA thesis, see \cite{O13} and \cite{J14}.
As we have seen above, when $\varGamma$ is a Lie groupoid over a manifold $M$, 
 $T\varGamma$ is a Lie groupoid over $TM$ and 
$T^*\varGamma$ is a Lie groupoid over $A^*$, where $A$ is the Lie algebroid of $\varGamma$.
Taking the direct sum of these Lie groupoids yields a Lie groupoid structure on the generalized tangent bundle ${\mathbb T} \varGamma = T\varGamma \oplus T^*\varGamma$ over $TM \oplus A^*$.
\begin{definition}
A Dirac structure, $L$, on a Lie groupoid, $\varGamma$, is {\it multiplicative} if it is defined by a subgroupoid of 
the Lie groupoid $T\varGamma \oplus T^*\varGamma \overrightarrow{\rightarrow} TM \oplus A^*$
over a sub-bundle ${\underline L} \subset TM \oplus A^*$.
A Lie groupoid equipped with a multiplicative Dirac structure is called a {\it Dirac groupoid}.
\end{definition}
\begin{example} If $\pi$ is a bivector on a
Lie groupoid, $\varGamma$, then $\pi^\sharp : T^*\varGamma \to T\varGamma$ is a morphism of groupoids over a morphism of vector bundles $A^* \to TM$
if and only if its graph is a subgroupoid of  $T\varGamma \oplus T^*\varGamma \overrightarrow{\rightarrow}  TM \oplus A^*$.
Therefore, in a Poisson groupoid, $(\varGamma,\pi)$, the graph of $\pi^\sharp$ defines a multiplicative Dirac structure.
The Poisson groupoids are those Dirac groupoids whose multiplicative Dirac structure is the graph of a morphism of Lie groupoids from $T^*\varGamma$ to $T\varGamma$.
\end{example} 
\begin{example} Presymplectic groupoids are Dirac groupoids whose multiplicative Dirac structure is the graph of a morphism of Lie groupoids from $T\varGamma$ to $T^*\varGamma$.
\end{example}
\begin{example} Dirac groups were defined and studied by 
Ortiz \cite{O08} and by Jotz \cite{J11} before the general concept of Dirac groupoids was introduced.
\end{example}
\subsection{The generalized tangent bundle of the Lie algebroid of a Lie groupoid}
Let $A$ be a Lie algebroid with base $M$. Then
$TA$ is a Lie algebroid with base $TM$,
and $T^*A$ is a Lie algebroid with base $A^*$. In fact, the direct sum {${\mathbb T} A = TA \oplus T^*A$} is a {Lie algebroid with base $TM \oplus A^*$} (see \cite{GU97} \cite{BZ09} \cite{JSX11}).
\begin{example}
A $2$-form, $\xi$, on a Lie algebroid, $A$, is morphic if and only if the graph of $\xi^\flat$ is a Lie subalgebroid of the Lie algebroid $TA \oplus T^*A \to TM \oplus A^*$.
\end{example}
\begin{proposition}
Let $\varGamma  \overrightarrow{\rightarrow} M$ be a Lie groupoid with Lie algebroid $A \to M$.
The {Lie algebroid of the Lie groupoid 
{$T\varGamma \oplus T^*\varGamma \overrightarrow{\rightarrow} TM \oplus A^*$}  is the Lie algebroid $TA \oplus T^*A \to TM \oplus A^*$}.
\end{proposition}
In fact, the Lie algebroid $TA \to TM$ is identified with the Lie algebroid of the tangent Lie groupoid, $T\varGamma  \overrightarrow{\rightarrow} TM$, and the Lie algebroid $T^*A \to A^*$ is identified with the Lie algebroid of the cotangent Lie groupoid, $T^*\varGamma \overrightarrow{\rightarrow} A^*$. These facts were proved in \cite{MX94} and this proposition follows. See also \cite{JSX11}.

\subsection{Dirac algebroids}
We shall now define morphic Dirac structures on Lie algebroids which are the infinitesimals of multiplicative Dirac structures on Lie groupoids.
\begin{definition}
A Dirac structure on a vector bundle, $E \to M$, is said to be {\emph{linear}} 
if it is defined by a vector sub-bundle of $TE \oplus T^*E \to TM \oplus E^*$. 
A Dirac structure on a Lie algebroid, $A \to M$, is said to be {\emph{morphic}} if it is defined by a Lie subalgebroid of $TA \oplus T^*A \to TM \oplus A^*$. 
A Lie algebroid equipped with a morphic Dirac structure is said to be a {\emph{Dirac algebroid}}. 
\end{definition}
\begin{remark} This terminology conflicts with the term used above for a Dirac bundle in a Courant algebroid. 
\end{remark}

The structure induced on the Lie algebroid, $A$, of a Lie groupoid, $\varGamma$, by a Dirac structure on $\varGamma$ is a Dirac structure on $A$. When the structure on $\varGamma$ is multiplicative, the structure on $A$ is infinitesimally multiplicative, i.e.,                                                                                                                                                                  \begin{theorem}
The infinitesimal of a Dirac groupoid is a Dirac algebroid. 
\end{theorem}

\begin{example}
If the Lie subalgebroid of $TA \oplus T^*A$ defining a Dirac algebroid structure on $A$ is the graph of a bivector, $\pi:  T^*A \to TA$, this bivector is a morphic Poisson bivector and therefore it defines a Lie bialgebroid $(A,A^*)$. More precisely, Lie bialgebroid structures on $A$ are in one-to-one correspondence  with Dirac algebroid structures on $A$ defined by the graph of a bivector.
\end{example}
\begin{example}
A necessary and sufficient condition for the graph of a $2$-form on $A$ to define a Dirac structure on $A$ is that it be closed and morphic. As a consequence, IM $2$-forms on $A$ are in one-to-one correspondence  with Dirac algebroid structures on $A$ defined by the graph of a $2$-form.
\end{example}
\begin{example}
The infinitesimal of a multiplicative Dirac structure on a Lie group, $G$, defines an ideal $\mathfrak{h}$ in the Lie algebra $\mathfrak{g}$ of $G$ such that $\mathfrak{g}/\mathfrak{h}$ and its dual form a Lie bialgebra.
\end{example}
\subsection{Excursus on ``Dirac-Lie algebroids''}
 With a view to applications in the mechanics of systems defined by implicit Euler-Lagrange or Hamilton equations, Grabowska and Grabowski introduced the concept of ``Dirac-Lie algebroids'' \cite{GG11}, which are vector bundles, $V \to M$, equipped with a linear Dirac structure on $V^*$.
This concept generalizes that of Lie algebroid because a Lie-algebroid structure on $V$ is equivalent to a linear Poisson structure on $V^*$ and its graph defines a linear Dirac structure on $V^*$. 

\begin{remark} The ``Dirac-Lie algebroids'' \cite{GG11}, that generalize the Lie algebroids, should not be confused with the ``Dirac algebroids'' \cite{O13} (for which the Dirac structure is on the vector bundle itself, not on its dual, and is assumed to be morphic) that generalize the Lie bialgebroids. The latter are the infinitesimals of Dirac groupoids, just as Lie bialgebroids are the infinitesimals of Poisson groupoids.  They are to be compared with the ``Dirac bialgebroids'' of \cite{J14}.
\end{remark}

\subsection{Generalized complex structures}
Recall that the {\it Nijenhuis torsion} of a $(1,1)$-tensor $N$ on a manifold, $M$, is the $(1,2)$-tensor, ${\mathcal T}_N$ such that
$$
{\mathcal T}_N(X,Y) = [NX,NY]- N[NX,Y] - N[X,NY] + N^2 [X,Y],
$$
for all vector fields $X$ and $Y$ on $M$, where $[~,~]$ is the Lie bracket of vector fields. 
A {\it complex structure} on a manifold, $M$, is an endomorphism of the tangent bundle whose square is $-{\mathrm{Id}}_{TM}$ and whose Nijenhuis torsion vanishes.

The {\it generalized complex structures} on a manifold, $M$, are defined similarly, replacing the tangent bundle by the generalized tangent bundle, ${\mathbb T} M$, and the Lie bracket of vector fields by the Dorfman bracket of sections of ${\mathbb T} M$. They are skew-symmetric endomorphisms of  ${\mathbb T} M$, whose square is $-{\rm {Id}}_{{\mathbb T} M}$, and whose Nijenhuis torsion vanishes. (The Nijenhuis torsion of a skew-symmetric endomorphism of ${\mathbb T} M$ of square $\pm {\rm {Id}}_{{\mathbb T} M}$ is again a tensor. See \cite{K11}.)

A generalized complex structure on $M$ is equivalent to a pair of transverse, complex conjugate Dirac structures in the complexified generalized tangent bundle of $M$.

A generalized complex structure on $M$ is of the form,
${\mathcal N} =  \begin{pmatrix}N & \pi^\sharp \\  \omega^\flat & - \, ^t  N \end{pmatrix}$,
where $N$ is a $(1,1)$-tensor, $\pi$ is a bivector and $\omega$ is a 2-form on $M$.

\subsection{Multiplicative generalized complex structures on Lie groupoids}
When $\G$ is a Lie groupoid over $M$, ${\mathbb T} \varGamma = T\varGamma \oplus T^*\varGamma$ is a Lie groupoid over $TM \oplus A^*$, where $A$ is the Lie algebroid of $\varGamma$. 
\begin{definition} 
A {generalized complex structure} on a Lie groupoid, $\varGamma$, is called \emph{multiplicative} if it is defined by a {Lie-groupoid automorphism} of ${\mathbb T} \varGamma$ over a vector-bundle automorphism of $TM \oplus A^*$.
\end{definition}
The concept of multiplicative generalized complex structure on a Lie groupoid was introduced and studied by Jotz, Sti\'enon and Xu \cite{JSX11}, who coined the name ``Glanon groupoids'' for them.
``Glanon groupoids'' comprise both symplectic groupoids and holomorphic groupoids.

\begin{example}[Glanon groups] A multiplicative generalized complex structure on a Lie group, $G$, is a groupoid automorphism of the Lie groupoid 
$TG \oplus T^*G \overrightarrow{\rightarrow} {\mathfrak g}^*$.
\end{example}   

\begin{example}[Symplectic groupoids] Let $\omega$ be a non-degenerate $2$-form on a Lie groupoid $\G$. Then the skew-symmetric endomorphism of $T\varGamma \oplus T^*\varGamma$, 
$ \begin{pmatrix}
 0 & - (\omega^\flat)^{-1} \\ \omega^\flat  & 0 \end{pmatrix}$, 
is a multiplicative generalized complex structure if and only if $(\varGamma, \omega)$ is a {symplectic groupoid}.
\end{example}

\begin{example}[Holomorphic groupoids] A holomorphic groupoid \cite{LSX09} is a Lie groupoid, $\varGamma$,  equipped with a vector-bundle endomorphism, ${{N}} : T\varGamma \to T\varGamma$, 
that satisfies the conditions: (i) ${{N}}^2 = - {\rm {Id}}_{T\varGamma}$, (ii) ${\mathcal{T}}_{{N}} = 0$, and (iii) ${{N}}$ is multiplicative, i.e., $N$ is a Lie groupoid automorphism of $T\G$ over a vector-bundle automorphism, ${\underline{{N}}} : TM \to TM$
(that satisfies ${\underline{{N}}}^2 = - {\rm{Id}}_{TM}$ and ${\mathcal T}_{\underline{{N}}} = 0$).

The {holomorphic Lie groupoid} structures, $N$, on $\varGamma$ coincide with those  multiplicative generalized complex structures which are skew-symmetric endomorphisms of $T\varGamma \oplus T^*\varGamma$ of the form,
${\mathcal{N}} = \begin{pmatrix}N &  0 \\ 0 &  - \, ^t \! N 
\end{pmatrix}$.
\end{example}

\begin{example}[Symplectic-Nijenhuis groupoids] Multiplicative generalized complex structures on $\G$ of the form,
$\N = \begin{pmatrix} N &  - (\omega^\flat)^{-1}  \\ \omega^\flat  & - \, ^t \!N \end{pmatrix} $,
where $N : T\varGamma \to T\varGamma$ and  $\omega$ is a non-degenerate 2-form, correspond to symplectic-Nijenhuis groupoids \cite{SX07} \cite{C11}.
\end{example}   

\begin{proposition} 
A multiplicative generalized complex structure on a Lie groupoid, $\G$, induces a Poisson groupoid structure on $\G$.
\end{proposition}

In fact a multiplicative generalized complex structure, $\N$, on $\G$ is necessarily of the form 
$\N = \begin{pmatrix} N & \pi^\sharp \\ \omega^\flat  & - \, ^t \!N \end{pmatrix}$, with $\pi^\sharp : T^*\varGamma \to T \varGamma$ a morphism of Lie groupoids, i.e. the associated bivector, $\pi$, is multiplicative, and $\pi$ is a Poisson bivector.

\subsection{The main results of Jotz-Sti\'enon-Xu}

What is the infinitesimal structure induced by a multiplicative generalized complex structure on $\varGamma$ 
on the Lie algebroid $TA \oplus T^*A \to TM \oplus A^*$, where $A$ is the Lie algebroid of $\varGamma$?

\begin{definition} A generalized complex structure, $J : TA \oplus T^*A \to TA \oplus T^*A$, on a Lie algebroid, $A$, is called {\it infinitesimally multiplicative} if it is a
{Lie-algebroid automorphism} of $TA \oplus T^*A$ over a vector-bundle automorphism, $\underline J : TM \oplus A^* \to TM \oplus A^*$ .
\end{definition}

\begin{example} 
Let $J$  be an infinitesimally multiplicative generalized complex structure  on a Lie algebroid, $A$, of the form $\begin{pmatrix} 0 & - (\xi^\flat)^{-1}  \\ \xi^\flat  & 0 \end{pmatrix}$. Then $\xi^\flat : TA \to T^*A$ defines an IM 2-form on $A$.
\end{example}

\begin{proposition} 
 If a Lie algebroid, $A$, is equipped with an infinitesimally multiplicative generalized complex structure, there is an associated Lie algebroid structure on $A^*$ such that  $(A,A^*)$ is a Lie bialgebroid.
\end{proposition}

\begin{theorem} 
Let $A$ be the Lie algebroid of a Lie groupoid, $\varGamma$. Then any multiplicative generalized complex structure on $\varGamma$ 
induces an infinitesimally multiplicative generalized complex structure on $A$.
\end{theorem}

While there are two constructions of the Lie algebroid structure of the pair $(A,A^*)$, they yield the same result.

\begin{theorem} 
Let $\varGamma$ be a Lie groupoid equipped with a multiplicative generalized complex structure, $\N$. Let $A$ be the Lie algebroid of $\varGamma$. When $\varGamma$ is equipped with the Poisson groupoid structure induced by $\N$, the Lie bialgebroid of $\varGamma$ is the Lie bialgebroid $(A,A^*)$ associated to the infinitesimally multiplicative generalized complex structure on $A$ induced by $\N$.
\end{theorem}
 
The following integration theorem subsumes diverse results such as the integration of Poisson manifolds into symplectic groupoids.
 \begin{theorem} 
If $\varGamma$ is a source-simply-connected Lie groupoid integrating a Lie algebroid, $A$, with an infinitesimally multiplicative generalized complex structure, then there is a unique multiplicative generalized complex structure on $\varGamma$ integrating that of $A$. 
\end{theorem}


\section*{Conclusion}

The re-formulation of the ``multiplicativity'' property of bivectors in terms of morphisms is a powerful tool
that has permitted identifying the adequate notion of {compatibiblity of a Poisson-like or a Nijenhuis-like structure with a group-like structure}, generalizing the multiplicativity property of the Poisson bivector on a Poisson-Lie group.
This is the program carried out in recent years whose results I have tried to outline here.
Some important developments have not been alluded to, such as contact structures, quasi-Poisson structures and multiplicative Manin pairs, double Lie groupoids, Spencer operators or multiplicativity-up-to-homotopy, and many other topics.
It remains to define and exploit the general concept of multiplicativity for multivector-valued forms of higher-degree that would encompass all the known cases of forms, multivectors and vector-valued $1$-forms.
(See a paper in preparation by Bursztyn and Thiago Drummond.)

\bigskip

\subsection*{Acknowledgments}

This text is an expanded version of a lecture delivered at the conference held at the Banach Center, ``Geometry of jets and fields'', in honor of Janusz Grabowski, May 10-16, 2015. I am grateful to the organizers for their invitation and for the pleasant stay at B{\k e}dlewo. I thank the referee for a careful reading of this paper and for a few, valuable suggestions.
My warmest thanks to Kirill Mackenzie for his help during the completion of this manuscript and to Camille Laurent-Gengoux for much information, including an unpublished result, that he communicated to me during one of our always instructive exchanges.

\end{document}